\documentclass[12pt]{amsart}
\usepackage{amssymb}
\usepackage{amsmath}
\usepackage{longtable}
\newcommand\g{{\mathfrak g}}
\newcommand\h{{\mathfrak h}}

\renewcommand{\b}{\mathfrak{b}}

\renewcommand{\a}{\mathfrak{a}}
\renewcommand{\c}{\mathfrak{c}}
\newcommand\m{\mathfrak m}
\newcommand\lfrak{\mathfrak l}
\newcommand\n{\mathfrak n}
\newcommand\q{\mathfrak q}
\newcommand\z{\mathfrak z}

\newcommand\p{\mathfrak p}

\renewcommand{\t}{\mathfrak{t}}

\newcommand\im{\operatorname{im}}
\newcommand\Spec{\operatorname{Spec}}
\newcommand\Quot{\operatorname{Quot}}

\newcommand\Aut{\operatorname{Aut}}

\newcommand\K{\mathbb K}
\newcommand\X{\mathfrak X}
\newcommand\Q{\mathbb Q}
\newcommand\R{\mathbb R}

\newcommand\N{\mathcal N}
\newcommand\Z{\mathbb Z}
\newcommand\A{\mathfrak A}

\newcommand\Rad{\operatorname{R}}

\newcommand\GL{\mathop{\rm GL}\nolimits}

\newcommand\SL{\mathop{\rm SL}\nolimits}

\newcommand\Span{\operatorname{Span}}

\newcommand{\Ad}{\mathop{\rm Ad}\nolimits}

\newcommand{\rank}{\mathop{\rm rk}\nolimits}

\newcommand{\defe}{\mathop{\rm def}\nolimits}

\renewcommand{\Ad}{\mathop{\rm Ad}\nolimits}
\newcommand{\Hom}{\operatorname{Hom}}

\newcommand\quo{/\!/}
\newtheorem{Thm}{Theorem}[section]
\newtheorem{Prop}[Thm]{Proposition}
\newtheorem{Cor}[Thm]{Corollary}
\newtheorem{Lem}[Thm]{Lemma}
\theoremstyle{definition}
\newtheorem{Ex}[Thm]{Example}
\newtheorem{defi}[Thm]{Definition}
\newtheorem{Rem}[Thm]{Remark}

\unitlength=1mm   \numberwithin{equation}{section}
\numberwithin{table}{section} \oddsidemargin=0cm
\author{Ivan V. Losev}
\title{Combinatorial invariants of algebraic Hamiltonian actions}
 \thanks{Supported
by RFBR grant  05-01-00988}
\begin{document}
\begin{abstract}
To any Hamiltonian action of a reductive algebraic group $G$ on a
smooth irreducible symplectic variety $X$ we associate certain
combinatorial invariants: Cartan space, Weyl group, weight and root
lattices.  For cotangent bundles these invariants essentially
coincide with those arising in the theory of equivarant embeddings.
Using our approach we establish some properties of the latter
invariants.
\end{abstract}
\maketitle \tableofcontents
\section{Introduction}
In this paper we study certain invariants associated with
Hamiltonian actions of reductive algebraic groups. All groups and
varieties are defined over an algebraically closed field $\K$ of
characteristic zero.

Let $G$ be a reductive algebraic group and $X$ a Hamiltonian
$G$-variety, i.e.,  a smooth symplectic variety equipped with a
Hamiltonian action of $G$ (see Section \ref{SECTION_prelim} for
generalities on Hamiltonian actions). To define our invariants we
need to fix a certain Levi subgroup $L\subset G$ and an $L$-stable
subvariety (a so-called {\it $L$-cross-section}) $X_L\subset X$ (see
Section \ref{SECTION_crossec} for definitions). To the pair
$(L,X_L)$ we associate the affine subspace $\a_{G,X}^{(X_L)}\subset
\z({\mathfrak l})$ (the {\it Cartan space} of $X$), the subgroup
$W_{G,X}^{(X_L)}\subset N_G(\a_{G,X}^{(X_L)})/Z_G(\a_{G,X}^{(X_L)})$
(the {\it Weyl group} of $X$), and the lattices
$\X_{G,X}^{(X_L)},\Lambda_{G,X}^{(X_L)}\subset \a_{G,X}^{(X_L)}$
(the {\it weight} and {\it root lattices}, respectively). In the
sequel, while speaking about combintorial invariants of $X$,  we
mean these data.

Our interest to these invariants comes from the theory of
equivariant embeddings of homogeneous spaces. Namely, suppose $G$ is
connected. Fix a Borel subgroup $B\subset G$ and a maximal torus $
T\subset B$. To an arbitrary irreducible $G$-variety $X_0$ one
associates its combintorial invariants: the weight lattice
$\X_{G,X_0}$, the Cartan space $\a_{G,X_0}$, the Weyl group
$W_{G,X_0}\subset \GL(\a_{G,X_0})$ and the root lattice
$\Lambda_{G,X_0}$. Let us give here a short description of these
invariants, precise definitions will be given in the beginning of
Section \ref{SECTION_cotangent}.

By definition, the weight lattice is nothing else but the set of all
weights $\lambda\in \X(T)$ such that there exists a
$B$-semiinvariant rational function on $X_0$ of weight $\lambda$.
The Cartan space $\a_{G,X_0}$ is a subspace of $\t^*$ spanned by
$\X_{G,X_0}$. The Weyl group describes the structure of the set of
so called {\it central} $G$-invariant valuations of $\K(X_0)$ (see
\cite{Brion},\cite{Knop4}). Similarly, the lattice $\Lambda_{G,X_0}$
introduced in \cite{Knop8} is, in a certain sense,  dual to the
group of {\it central} $G$-equivariant automorphisms of $\K(X_0)$.

There is a special class of $G$-varieties for that the role of the
above invariants is particularly important. Namely, a normal
$G$-variety $X_0$ is called {\it spherical} if $B$ has an open orbit
on $X_0$. In this case all $G$-valuations and all $G$-automorphisms
are central. The Cartan space and the Weyl group play an essential
role in the classification theory of spherical varieties
generalizing that of toric varieties, see, for example, \cite{Knop5}
for details.

Two sets of invariants are closely related. Namely, to a smooth
$G$-variety $X_0$ one associates its cotangent bundle $X$. The
latter is equipped with the natural symplectic form and the action
of $G$. This action is Hamiltonian. For $L$ one takes
$Z_G(\a_{G,X_0})$. Further, Vinberg, \cite{Vinberg}, essentially
noticed that there is a distinguished choice of the subvariety
$X_L$. It turns out that $\a_{G,X_0}=\a_{G,X}^{(X_L)},
\X_{G,X_0}=\X_{G,X}^{(X_L)}, W_{G,X_0}=W_{G,X}^{(X_L)},
\Lambda_{G,X_0}=\Lambda_{G,X}^{(X_L)}$ provided $X_0$ is
quasiaffine. The first two equalities follow easily from the
construction of Sections 2,3 in \cite{Knop3}. The equality of the
Weyl groups is the main result of \cite{Knop3}, the proof is rather
complex. Finally, the equality for root lattices was implicitly
proved in \cite{Knop8}.

Now we describe the content of the paper. Section
\ref{SECTION_prelim} contains preliminaries on Hamiltonian
varieties. In Section \ref{SECTION_crossec} we review the theory of
cross-sections of Hamiltonian actions tracing back to Guillemin and
Sternberg, \cite{GS}. An example of a cross-section is provided by a
subvariety $X_L$ mentioned in the beginning of the present section.

Section \ref{SECTION_morph} is devoted to the notion of a
Hamiltonian morphism between two Hamiltonian varieties. Using the
notion of a cross-section   we define  {\it central} automorphisms
of a Hamiltonian variety. Then we study some of their properties.

In Section \ref{SECTION_combinv} we define the combintorial
invariants of a Hamiltonian variety  and study their basic
properties, Lemmas \ref{Lem:2.7.9}-\ref{Lem:2.7.6}.

In Section \ref{SECTION_cotangent} we compare the combintorial
invariants of a $G$-variety $X_0$ with those defined for the
Hamiltonian action on its cotangent bundle. First of all, we recall
the definitions of the combinatorial invariants for $X_0$. Then we
state the main result of this section, the comparision theorem
\ref{Thm:6.4}. Next we define the distinguished cross-sections
$\Sigma_M$ of $X$ for certain Levi subgroups $M\subset G$ and study
their properties (Proposition \ref{Prop:3.1.3}). Finally we prove
Theorem \ref{Thm:6.4}. The proof is contained in Propositions
\ref{Prop:3.1.4}, \ref{Prop:3.2.3}. None of this results is actually
new. Note however that our proof of the equality of the root
lattices (Proposition \ref{Prop:3.2.3}) is much more elementary than
that given in \cite{Knop8} (Theorems 4.1, 7.8). The proof in
\cite{Knop8} was valid for arbitrary (i.e., not necessarily
quasiaffine) $G$-varieties but used the theory of group schemes in
an essential way.

In Section \ref{SECTION_properties} we prove some results about
$W_{G,X_0}, \X_{G,X_0},\Lambda_{G,X_0}$ based on the comparision
theorem and results of Section \ref{SECTION_combinv}. The most
important results are Propositions \ref{Prop:3.4.2},\ref{Prop:3.4.4}
and Theorem \ref{Thm:3.0.3}. Proposition \ref{Prop:3.4.2} allows one
to compare the Weyl groups $W_{G,X_0}$ and $W_{M,X_0/\Rad_u(Q)}$,
where $Q$ is a parabolic subgroup of $G$ containing $B$ and subject
to some additional conditions, and $M$ is the Levi subgroup of $Q$
containing $T$. This result plays an important role, for example, in
the proof of the Knop conjecture, see \cite{Knop_conj}. Theorem
\ref{Thm:3.0.3} reduces the computation of Weyl groups and weight
lattices to the case when $\a_{G,X_0}=\t$. Proposition
\ref{Prop:3.4.4} continues the reduction procedure. In the case when
$\a_{G,X_0}=\t$ it reduces the computation of
$W_{G,X_0},\Lambda_{G,X_0}$ to the case when $G$ is simple.

Results and constructions of this paper not connected with root
lattices constitute Sections 3,4 of the preprint \cite{Weyl}. They
were used there to construct algorithms computing Weyl groups of
$G$-varieties. The theory of root lattices developed here will be
used in a forthcoming paper to establish an algorithm computing {\it
weight} lattices of $G$-varieties.

\section{Notation and conventions}\label{SECTION_notation}
If an algebraic group is denoted by a capital Latin letter, then we
denote its Lie algebra by the corresponding small German letter.
\begin{longtable}{p{5.5cm} p{10cm}}
$\langle\alpha,v\rangle$& the pairing between elements $\alpha,v$ of
dual vector spaces.
\\ $A^\times$& the group of invertible elements of an associative
algebra $A$.
\\ $Aut^G(X)$& the group of all $G$-equivariant automorphisms of a
$G$-variety $X$.
\\ $\defe_G(X)$& the defect of an irreducible Hamiltonian
$G$-variety $X$.
\\ $(G,G)$ (resp., $[\g,\g]$)& the commutant of an algebraic group
$G$ (resp., of a Lie algebra $\g$)
\\ $G^{\circ}$& the connected component of unit of an algebraic
group $G$.
\\ $G*_HY$& the homogeneous bundle over $G/H$ with a fiber $Y$.
\\ $[g,y]$& the equivalence class of $(g,y)$ in $G*_HY$.
\\ $G_x$& the stabilizer of a point $x\in X$ under an action
$G:X$.
\\ $m_G(X)$& the maximal dimension of an orbit for an action of an
algebraic group $G$ on a variety $X$.
\\ $N_G(H)$, (resp. $N_G(\h),\n_\g(\h)$)& the normalizer of an
algebraic subgroup $H$ in an algebraic group $G$ (resp. of a
subalgebra $\h\subset \g$  in an algebraic group $G$, of a
subalgebra $\h\subset \g$ in a Lie algebra $\g$).
\\ $N_G(H,Y)$& the stabilizer of a subset $Y\subset X$ under the
action $N_G(H):Y$, where $X$ is a set acted on by $G$.\\
$\Rad_u(G)$& the unipotent radical of an algebraic group
$G$\\$v(f)$& the skew-gradient of a rational function on a
symplectic
variety $X$.\\
$V^{(G)}_\lambda$&$=\{v\in V| gv=\lambda(g)v, \forall g\in G\}$,
where $G$ is an algebraic group, $V$ is a $G$-module, $\lambda\in
\X(G)$.\\
 $W(\g)$& the Weyl group of a reductive Lie algebra $\g$.
\\$\X(G)$& the character group of an algebraic group $G$.\\
$X^G$& the
fixed-point set
for an action of $G$ on $X$.\\
$X\quo G$& the categorical quotient for an action of a reductive
group on
an affine variety $X$\\
$X^{reg}$& the subset of an algebraic variety $X$ consisting of all smooth points.\\
 $\#X$& the cardinality of a
set $X$.
\\ $\z(\g)$& the center of a Lie algebra $\g$.
\\   $Z_G(H)$, (resp. $Z_G(\h),\z_\g(\h)$)& the centralizer of an
algebraic subgroup $H$ in an algebraic group $G$ (resp. of a
subalgebra $\h\subset \g$  in an algebraic group $G$, of a
subalgebra $\h\subset \g$ in a Lie algebra $\g$).\\
$\Lambda(\g)$& the root lattice of a reductive Lie algebra $\g$.\\
$\mu_{G,X}$& the moment map of a Hamiltonian $G$-variety $X$.\\
$\xi_s$& the semisimple part of an element $\xi\in \g$.\\
$\pi_{G,X}$& the quotient morphism $X\rightarrow X\quo G$ for the
action $G:X$.
\end{longtable}
\section{Preliminaries}\label{SECTION_prelim}
Here we recall some basic definitions and results concerning
Hamiltonian actions.

In this section $X$ is a smooth algebraic variety equipped with a
regular symplectic form $\omega$ and $G$ is a reductive algebraic
group acting on $X$ by symplectomorphisms.

 Let $U$ be an open subset of $X$ and $f$ a regular
function on $U$.
 The {\it skew-gradient} $v(f)$ of $f$ is a regular vector field on $U$
 given by the equality
\begin{equation}\label{eq:2.1:1}
\omega_x(v(f),\eta)=\langle d_xf, \eta\rangle, x\in U, \eta\in T_xX.
\end{equation}
 For $f,g\in \K[U]$ one can define their Poisson bracket
$\{f,g\}\in \K[U]$ by
\begin{equation}\label{eq:2.1.2}
\{f,g\}=\omega(v(f),v(g)).
\end{equation}
Clearly, $\{f,g\}=L_{v(f)}g$, where $L$ denotes the Lie derivative.

To any element $\xi\in\g$ one associates the velocity vector field
$\xi_*$. Suppose  there is a linear map $\g\rightarrow \K[X],
\xi\mapsto H_\xi,$ satisfying the following two conditions:
\begin{itemize}
\item[(H1)] The map $\xi\mapsto H_\xi$ is $G$-equivariant.
\item[(H2)] $v(H_\xi)=\xi_*$.
\end{itemize}

\begin{defi}\label{Def:2.1.1}
The action  $G:X$ equipped with a linear map $\xi\mapsto H_\xi$
satisfying (H1),(H2) is said to be {\it Hamiltonian} and $X$ is
called a Hamiltonian $G$-variety. 
\end{defi}

For a Hamiltonian action $G:X$ we define the morphism
$\mu_{G,X}:X\rightarrow \g^*$ by the formula
\begin{equation}\label{eq:2.1:3}
\langle \mu_{G,X}(x),\xi\rangle= H_{\xi}(x),\xi\in\g,x\in X.
\end{equation}
This morphism is called the {\it moment map} of the Hamiltonian
$G$-variety $X$.

Conditions (H1),(H2) are equivalent, respectively, to

\begin{itemize}
\item[(M1)] $\mu_{G,X}$ is a $G$-morphism.
\item[(M2)] $\langle d_x\mu_{G,X}(v),\xi\rangle=\omega_x(\xi_*x,v),$ for all $ x\in X,v\in
T_xX,\xi\in\g$.
\end{itemize}

Any two maps $\mu_{G,X}:X\rightarrow \g^*$ satisfying  conditions
(M1),(M2)  differ by an element of $\g^{*G}$. Moreover,
$H_{[\xi,\eta]}=\{H_\xi, H_\eta\} =\omega(\xi_*,\eta_*)$ (see, for
example, \cite{Vinberg}).

The following example plays a crucial role in this paper.

\begin{Ex}[Cotangent bundles]\label{Ex:2.1.4} Let $Y$ be a smooth $G$-variety. Let $X=T^*Y$ be the cotangent
bundle of $Y$. $X$ is a symplectic algebraic variety (the symplectic
form is written down, for example, in \cite{Vinberg}). The action of
$G$ on $X$ is Hamiltonian. The moment map is given by
$\langle\mu_{G,X}((y,\alpha)), \xi\rangle=\langle \alpha,
\xi_{*}y\rangle$. Here $y\in Y, \alpha\in T^*_yY,\xi\in\g$.
\end{Ex}

\begin{Ex}[Symplectic vector spaces]\label{Ex:2.1.5}  Let $V$ be a vector space equipped with a non-degenerate
skew-symmetric bilinear form $\omega$. Then $V$ is a symplectic
variety. Let $G$ be a reductive group acting on $V$ by linear
symplectomorphisms. Then the action $G:V$ is Hamiltonian. The moment
map $\mu_{G,V}$ is given by  $\langle\mu_{G,V}(v), \xi\rangle =
\frac{1}{2}\omega(\xi v,v), \xi\in\g, v\in V$.
\end{Ex}

One fixes a $G$-invariant symmetric form $(\cdot,\cdot)$ on $\g$
such that the restriction of  $(\cdot,\cdot)$ to $\t(\Q)$ is a
scalar product. Such a form is necessarily nondegenerate so we may
identify $\g$ and $\g^*$ with respect to it.   In particular, we may
consider the moment map as a morphism from $X$ to $\g$.

\begin{Rem}[the restriction to a subgroup]\label{Ex:2.1.9} Let $H$ be a reductive subgroup
of $G$ and $X$ a Hamiltonian $G$-variety. Then $X$ is a Hamiltonian
$H$-variety with the moment map $\mu_{H,X}=p\circ\mu_{G,X}$. Here
$p$ denotes the restriction map $p:\g^*\rightarrow \h^*$.
\end{Rem}

Finally, we recall the definition of the defect of a Hamiltonian
action. The proof of the next lemma is contained, for instance, in
\cite{Vinberg}, Chapter 2, Section 2. We put
$\psi_{G,X}:=\pi_{G,\g}\circ\mu_{G,X}$.

\begin{Lem}\label{Lem:2.1.10}
Suppose $X$ is irreducible. The following numbers are equal:
\begin{enumerate}
\item $\dim\ker\omega_x|_{\g_*x}$ for $x\in X$ in general position.
\item $m_G(X)-m_G(\overline{\im\mu_{G,X}})$.
\item $\dim\overline{\im\psi_{G,X}}$.
\end{enumerate}
These equal numbers are called the {\rm defect} of $X$.
\end{Lem}

We denote the defect of $X$ by $\defe_G(X)$.

\section{Local cross-sections}\label{SECTION_crossec}
In this subsection $X$ is a quasiprojective Hamiltonian
$G$-variety, $\omega$ is the symplectic form on $X$.

The local cross-section theorem reduces the study of a Hamiltonian
$G$-variety in an etale neighborhood of a point $x\in X$ to the case
when  $\mu_{G,X}(x)_s\in\z(\g)$.

Let $L$ be a Levi subgroup of $G$ and ${\mathfrak l}$ the corresponding Lie
algebra. Put ${\mathfrak l}^{pr}=\{\xi\in{\mathfrak l}| \z_\g(\xi_s)\subset {\mathfrak l}\}$. It is a
complement to the  set of zeroes of some element from
$K[{\mathfrak l}]^{N_G({\mathfrak l})}$, see, for example, \cite{alg_hamil}, Subsection
5.1.

Put $Y=\mu_{G,X}^{-1}({\mathfrak l}^{pr})$. Clearly, $Y$ is an $N_G(L)$-stable
locally-closed subvariety of $X$.
\begin{Prop}[\cite{Knop2}, Theorem 5.4 and \cite{alg_hamil},
Corollary 5.3 and Propositions 5.2, 5.4]\label{Prop:2.3.1} In the
preceding notation
\begin{enumerate}
\item
 $T_yX={\mathfrak l}^\perp_*y\oplus T_yY$
is a skew-orthogonal direct sum  for any $y\in Y$. In particular,
$Y$ is a smooth subvariety of $X$ and the restriction of $\omega$ to
$Y$ is nondegenerate. Thus $Y$ is equipped with the symplectic
structure.
\item The action $N_G(L):Y$ is Hamiltonian with the
moment map $\mu_{G,X}|_Y$. \item The natural morphism
$G*_{N_G(L)}Y\rightarrow X$ is etale. Its image is saturated.
\end{enumerate}
\end{Prop}  Recall that a subset $Z^0$ of a $G$-variety
$Z$ is called {\it saturated} if there exist a $G$-invariant
morphism $\varphi:Z\rightarrow Z_0$ and a subset $Z_0^0\subset Z_0$
such that $Z^0=\varphi^{-1}(Z_0^0)$.

\begin{defi}\label{Def:2.3.2}
An irreducible (=connected) component of $\mu_{G,X}^{-1}({\mathfrak l}^{pr})$
equipped with the structure of a Hamiltonian $L$-variety obtained by
the restriction of the Hamiltonian structure from
$\mu_{G,X}^{-1}({\mathfrak l}^{pr})$ is called an {\it $L$-cross-section} of
$X$.
\end{defi}

Under some special choice of $L$, Proposition~\ref{Prop:2.3.1} can
be strengthened.

\begin{defi}\label{Def:2.3.3}
The centralizer in $G^\circ$ of $\mu_{G,X}(x)_s$ for $x\in X$ in
general position is called the {\it principal centralizer} of $X$.
\end{defi}

The principal centralizer is defined uniquely up to
$G^\circ$-conjugacy (see~\cite{alg_hamil}, Subsection 5.2).

\begin{Prop}[\cite{alg_hamil}, Proposition 5.7]\label{Prop:2.3.4}
Let $L$ be the principal centralizer of $X$. The natural morphism
$G*_{N_G(L)}\mu_{G,X}^{-1}(\lfrak^{pr})\rightarrow X$ is an open
embedding and the group $N_G(L)$ acts transitively on the set of
irreducible components of $\mu_{G,X}^{-1}(\lfrak^{pr})$.
\end{Prop}

Till the end of the subsection $L$ is the principal centralizer and
$X_L$ is an $L$-cross-section of $X$. In the sequel we denote the
open saturated subset $G\mu_{G,X}^{-1}(\lfrak^{pr})\subset X$ by
$X^{pr}\index{xpr@$X^{pr}$}$.

\begin{Lem}\label{Lem:2.3.5}
The following conditions are equivalent:
\begin{enumerate}
\item $m_G(X)=\dim G$.
\item $\defe_G(X)=\rank G$.
\item $\overline{\im \mu_{G,X}}=\g$.
\item $L$ is a maximal torus of $G$ and $m_L(X_L)=\defe_L(X_L)=\rank G$.
\end{enumerate}
\end{Lem}
\begin{proof}
The equivalence of the first three conditions is quite standard,
see, for example, \cite{alg_hamil}, Proposition 3.7, Corollary 3.10.
It remains to show that (4) is equivalent to (1)-(3). It follows
from (3) that $L$ is a maximal torus. On the other hand, if $L$ is a
maximal torus, then (1)$\Leftrightarrow$(4), thanks to  Proposition
\ref{Prop:2.3.1}.
\end{proof}

Now we discuss the behavior of the principal centralizer and a
corresponding cross-section under some operations with Hamiltonian
varieties.

\begin{Lem}\label{Lem:2.3.7}
Let $X$ be a Hamiltonian $G$-variety, $L$ its principal centralizer
and $X_L$ its $L$-cross-section.
\begin{enumerate}
\item Let $M$ be a Levi subgroup in $G$ containing $L$. Then there
is  a unique $M$-cross-section $X_M$ of $X$ containing $X_L$, $L$ is
the principal centralizer  and $X_L$ is an $L$-cross-section of
$X_M$. Conversely, any $L$-cross-section of $X_M$ is also an
$L$-cross-section of $X$.
\item Let $X'$ be a Hamiltonian $G$-variety and $\varphi:X\rightarrow X'$
a generically finite dominant rational mapping satisfying
$\mu_{G,X'}\circ\varphi=\mu_{G,X}$. Then $L$ is the principal
centralizer of $X'$, the rational mapping $\varphi$ is defined in
$X_L$ and there exists a unique $L$-cross-section $X_L'$ of $X'$
such that $\overline{\varphi(X_L)}=\overline{X_L'}$.
\item Suppose $G$ is connected. Let $G^0\subset G$ be a connected subgroup  containing
$(G,G)$. Then $L^0:=G^0\cap L$ is the principal centralizer and
$X_L$ is an $L^0$-cross-section of the Hamiltonian $G^0$-variety
$X$.
\item Suppose $G$ is connected and $X$ satisfies the equivalent conditions
of Lemma \ref{Lem:2.3.5}. Let $G_1,\ldots,G_k$ be all simple normal
subgroups of $G$. Then the Hamiltonian $G_i$-variety $X$ also
satisfies the equivalent conditions of Lemma \ref{Lem:2.3.5}. Put
$G^{(i)}=Z(G)^\circ \prod_{j\neq i}G_j$. There exists a unique
$L_i:=L\cap G_i$-cross-section $X_{L_i}$ of $X$ containing
$G^{(i)}X_L$. Moreover, $G^{(i)}X_L$ is an open subvariety in
$X_{L_i}$.
\end{enumerate}
\end{Lem}
\begin{proof}
{\it Assertion 1.} Since $\lfrak^{pr}\subset\m^{pr}$, we see that $X_M$
exists. It is unique because $\mu_{G,X}^{-1}(\m^{pr})$ is smooth.
Note that $\z_{\m}(\mu_{G,X}(x)_s)=\lfrak$ for all $x\in X_L$. It
follows that $L$ contains the principal centralizer of $X_M$. By
assertion 3 of Proposition \ref{Prop:2.3.1}, $G\im\mu_{M,X_M}$ is
dense in $\overline{\im\mu_{G,X}}$. Therefore the principal
centralizers of $X_M$ and $X$ are conjugate in $G$. Thence $L$ is
the principal centralizer of $X_M$. Since $\lfrak^{pr}\subset\m^{pr}$,
we see that $X_L$ is an  $L$-cross-section of $X_M$. Since all
$L$-cross-sections of $X_M$ are $N_M(L)$-conjugate, we obtain the
remaining claim.

{\it Assertion 2.} Note that
$\overline{\im\mu_{G,X}}=\overline{\im\mu_{G,X'}}$. Thus $L$ is the
principal centralizer of $X'$. Since $GX_L$ is dense in $X$,
$\varphi$ is defined in $X_L$. Let $X'_L$ be a unique
$L$-cross-section of $X$ containing a point from $\varphi|_{X_L}$.
Since $\mu_{G,X'}^{-1}(\lfrak^{pr})$ is smooth,
$\overline{\varphi(X_L)}\subset \overline{X_L'}$. To prove the
equality it is enough to compare dimensions using assertion 3 of
Proposition \ref{Prop:2.3.1}.

{\it Assertion 3.} Straightforward. 

{\it Assertion 4.} We may change $G$ by some covering and assume
that $G=Z(G)^\circ \times G_1\times\ldots\times G_k$. It can be
easily seen that $L_i$ is the principal centralizer of the
Hamiltonian $G_i$-variety $X$. Note that $\mu_{G_i,X}$ is a
$G^{(i)}$-invariant morphism and $X_L\subset
\mu_{G_i,X}^{-1}(\lfrak_i^{pr})$. Thus $G^{(i)}X_L\subset
\mu_{G_i,X}^{-1}(\lfrak_i^{pr})$. Since $\mu_{G_i,X}^{-1}(\lfrak_i^{pr})$ is
smooth and $G^{(i)},X_L$ are irreducible, there exists a unique
$L_i$-cross-section $X_{L_i}$ of $X$ containing $G^{(i)}X_L$. We
recall that the natural morphism $G*_LX_L\rightarrow X$ is an open
embedding. Thus $G^{(i)}X_L$ is a locally-closed subvariety in $X$
and thence in $X_L$. Moreover, $\dim G^{(i)} X_L=\dim G^{(i)}-\rank
G^{(i)}+\dim X_L= \dim X-\dim G_i+\rank G_i=\dim X_{L_i}$.
\end{proof}

The following lemma is straightforward.

\begin{Lem}\label{Lem:2.3.8}
Let $X$ be a Hamiltonian $G$-variety, $M\subset G$ a Levi subgroup,
$X_M$ an $M$-cross-section of $X$. Suppose that there is an action
$\K^\times:X$ and an integer $k$ such that
$\mu_{G,X}(tx)=t^k\mu_{G,X}(x), x\in X$. Then $X_M$ is
$\K^\times$-stable and $\mu_{M,X_M}(tx)=t^k\mu_{M,X_M}(x), x\in
X_M$.
\end{Lem}

\section{Hamiltonian morphisms}\label{SECTION_morph}
\begin{defi}\label{Def:3.1}
Let $X,Y$ be Hamiltonian $G$-varieties. An  etale $G$-equivariant
morphism $\varphi: X\rightarrow Y$ is called {\it Hamiltonian} if it
satisfies conditions 1,2  below:
\begin{enumerate}
\item $\varphi^*\omega^Y=\omega^X$. Here $\omega^X,\omega^Y$ are
symplectic forms on $X$ and $Y$.
\item $\mu_{G,Y}\circ\varphi=\mu_{G,X}$.
\end{enumerate}
\end{defi}

\begin{Ex}\label{Ex:3.2}
Let $X_0$ be a smooth $G$-variety and $\varphi$ its automorphism.
Let $\varphi_*$ denote the natural symplectomorphism of $T^*X_0$
lifting $\varphi$. Note that the map $\varphi\mapsto \varphi_*$ is a
monomorphism. If $\varphi$ is a $G$-equivariant, then $\varphi_*$ is
a Hamiltonian automorphism of $X$.
\end{Ex}

In Section \ref{SECTION_cotangent} we use the the following
description of automorphisms from the previous example. Note that
there is the natural action $\K^\times:X:=T^*X_0$ by the fiberwise
multiplication. This action satisfies the assumptions of Lemma
\ref{Lem:2.3.8} for $k=1$.

Now suppose that $X_0$ is quasiaffine and $\phi:X\dashrightarrow X$
is a $G\times \K^\times$-equivariant rational mapping preserving the
symplectic form such that $\phi^*(\K[X])=\K[X]$. Then $\phi^*$
induces an automorphism of $\K[X_0]=\K[X]^{\K^\times}$.

\begin{Lem}\label{Lem:3.3}
Let $X,X_0,\varphi$ be such as above. Suppose that the restriction
of $\varphi^*$ to $\K[X_0]$ is induced by some regular automorphism
$\psi$ of $X_0$. Then $\varphi=\psi_*$.
\end{Lem}
\begin{proof}
Replacing $\varphi$ with $\varphi\psi_*^{-1}$ we may assume that
$\psi=id$. Note that $\varphi^*$ preseves the Poisson bracket on
$X$. The $\K[X_0]$-module $\K[X]^{(\K^\times)}_{-1}$ is identified
with the module of vector fileds on $X_0$. For $v\in
\K[X_0]^{(\K^\times)}_{-1}, g\in \K[X_0]$ we have $\{v,f\}=L_vf$. It
follows that $\varphi^*$ is the identity on
$\K[X_0]^{(\K^\times)}_{-1}$. Since $X_0$ is quasiaffine, there is
$f\in \K[X_0]$ such that $\{x\in X_0| f(x)\neq 0\}$ is affine. It
follows that $\K[X_0]$ and $\K[X_0]^{(\K^\times)}_{-1}$ generate
$\K(X)$. Therefore $\varphi=id$.
\end{proof}

Now we consider a special class of Hamiltonian automorphisms. Let
$L$ be the principal centralizer and $X_L$ an $L$-cross-section of
$X$.

\begin{defi}\label{Def:3.8}
A Hamiltonian automorphism $\varphi$ of $X$ is said to be {\it
central}, if $\varphi(X_L)=X_L$ and there exists $t_\varphi\in
Z(L)^\circ$ such that $\varphi(x)=t_\varphi x$ for all $x\in X_L$.
Central automorphisms form a group denoted by $\A_{G,X}^{(\cdot)}$.
\end{defi}

Remark \ref{Rem:2.5.7} in the next section shows that the condition
of being central does not depend on the choice of the principal
centralizer and a cross-section.



The following lemma partially justifies   the term "central".

\begin{Lem}\label{Lem:3.4}
Let $\psi$ be a Hamiltonian automorphism of $X$ fixing $X_L$. Then
$\psi$ commutes with any central Hamiltonian automorphism of $X$.
\end{Lem}
\begin{proof}
Let $\varphi$ be a central Hamiltonian automorphism of $X$. Since
$\psi,\varphi$ are $G$-equivariant and $GX_L$ is dense in $X$, it is
enough to show that
$\psi|_{X_L}\varphi|_{X_L}=\varphi|_{X_L}\psi|_{X_L}$. It remains to
note that the restriction of $\psi$ to $X_L$ is
$Z(L)^\circ$-equivariant.
\end{proof}

Let $A_{G,X}^{(X_L)}$ denote the quotient of $Z(L)^\circ$ by the
inefficiency kernel of the action $Z(L)^\circ:X_L$. The action
$N_G(L,X_L)$ on $Z(L)^\circ$ descends to the action
$N_G(L,X_L):A_{G,X}^{(X_L)}$.

The restriction of a central automorphism to $X_L$ determines the
homomorphism $\A_{G,X}^{(\cdot)}\rightarrow A_{G,X}^{(X_L)}$. Since
$GX_L$ is dense in $X$, this homomorphism is injective. Let us
denote its image by $\A_{G,X}^{(X_L)}$. Since the restriction of any
$G$-equivariant automorphism to $X_L$ is $N_G(L,X_L)$-equivariant,
we see that $\A_{G,X}^{(X_L)}\subset
(A_{G,X}^{(X_L)})^{N_G(L,X_L)}$.

Now we give a more convenient description of the subgroup
$\A_{G,X}^{(X_L)}\subset (A_{G,X}^{(X_L)})^{N_G(L,X_L)}$. Any
element $l$ of the last group defines a $G$-equivariant automorphism
$l_*$ of $X^{pr}=G*_{N_G(L,X_L)}X_L$ by $l_*[g,y]=[g,ly]$.

\begin{Lem}\label{Lem:3.5} Let $l\in (A_{G,X}^{(X_L)})^{N_G(L,X_L)}$. Then
\begin{enumerate}
\item The automorphism $l_*$ of $X^{pr}$ is Hamiltonian and central. So $(A_{G,X}^{(X_L)})^{N_G(L,X_L)}=\A_{G,X^{pr}}^{(X_L)}$.
\item $l\in \A_{G,X}^{(X_L)}$ iff $l_*$ can be extended to a morphism $X\rightarrow X$.
\end{enumerate}
\end{Lem}
\begin{proof}
The equality $\mu_{G,X^{pr}}\circ l_*=\mu_{G,X^{pr}}$ stems from
$\mu_{G,X}|_{X_L}=\mu_{L,X_L}$. To complete the proof of assertion 1
it remains to verify that $l_*$ preserves $\omega$. By assertion 1
of Proposition \ref{Prop:2.3.1}, for any $y\in X_L$ there is the
skew-orthogonal direct sum decomposition $T_yX=T_yX_L\oplus
\lfrak_*^{\perp}y$. Since $X_L$ is $l_*$-stable and $l_*$ is
$G$-equivariant, we have
 $d_y(l_*)T_yX_L=T_{ly}X_L, d_y(l_*)\lfrak_*^{\perp}y=\lfrak_*^\perp(ly)$. Since
 $y\mapsto ly$ is a symplectomorphism of $X_L$, we see that $d_y(l_*)$ maps
$\omega_{y}|_{T_yX_L}$ to $\omega_{ly}|_{T_{ly}X_L}$. On the other
hand, for $\xi,\eta\in \lfrak^{\perp}$ the equality
$\omega_y(\xi_*y,\eta_*y)=(\mu_{L,X_L}(y), [\xi,\eta])$ holds. Since
$y\mapsto ly$ is a Hamiltonian automorphism of $X_L$, it follows
that
$d_y(l_*)\omega_{y}|_{\lfrak^\perp_*y}=\omega_{ly}|_{\lfrak^\perp_*(ly)}$.

Proceed to assertion 2. If $l\in \A_{G,X}^{(X_L)}$, then there is
nothing to prove. Now let $l\in (A_{G,X}^{(X_L)})^{N_G(L,X_L)}$ be
such that $l_*$ can be extended to a morphism $\varphi:X\rightarrow
X$. By assertion 1, $\varphi^*\omega=\omega$. Since $\omega$ is
non-degenerate, we see that $\varphi$ is etale. Taking into account
that $\varphi$ is birational, we see that it is an open embedding.
But any open embedding $X\rightarrow X$ is an isomorphism (see, for
example, \cite{Iitaka}, Proposition 4).
\end{proof}

\begin{Cor}\label{Cor:3.6}
If $X$ is quasiaffine, then $\A_{G,X}^{(X_L)}$ is a closed subgroup
of $A_{G,X}^{(X_L)}$.
\end{Cor}
\begin{proof}
The subgroup
$\A_{G,X^{pr}}^{(X_L)}=(A_{G,X}^{(X_L)})^{N_G(L,X_L)}\subset
A_{G,X}^{(X_L)}$ is closed. By assertion 1 of Lemma \ref{Lem:3.5},
the corollary stems from the following lemma.\end{proof}
\begin{Lem}\label{Lem:3.7}
Let $G$ be an algebraic group acting on an open subset $X^0$ of a
quasiaffine variety $X$. Then the subgroup of $G$ consisting of all
elements, that can be extended to an automorphism of $X$, is closed.
\end{Lem}
\begin{proof} $\K[X^0]$ is a rational $G$-algebra (see \cite{VP}, Lemma 1.4)
and $\K[X]$ is a $G$-submodule of $\K[X^0]$. The set of all elements
$g\in G$ such that $g\K[X]\subset \K[X]$ is given by the condition
that some matrix coefficients vanish. Thence this set is closed.
Thus $\{g\in G| g\K[X]=\K[X]| g\K[X]\subset \K[X],
g^{-1}\K[X]\subset\K[X]\}$ is a closed subgroup of $G$. So we may
assume that $g\K[X]=\K[X]$ for all $g\in G$.

There is a finitely generated subgalgebra $A\subset \K[X]$ such that
the corresponding map $X\rightarrow \overline{X}:=\Spec(A)$ is an
embedding. Since $\K[X]$ is rational, we may change $A$ by a
subalgebra generated by $GA$ and assume that $A$ is $G$-stable. An
element $g\in G$ defines an automorphism of $X$ iff
$g(\overline{X}\setminus X)=\overline{X}\setminus X$. But the
stabilizer of a subvariety is  closed.
\end{proof}
\section{Combinatorial invariants}\label{SECTION_combinv}
Throughout this subsection $X$ is an irreducible quasiprojective
Hamiltonian $G$-\!\! variety. We denote by $L$ the principal
centralizer of $X$.

Choose an $L$-cross-section  $X_L$ of $X$. Note that the triple
$(L,X_L,N_G(L,X_L))$ is defined uniquely up to $G$-conjugacy. Denote
by $T_0$ the inefficiency kernel for the action $Z(L)^\circ:X_L$ and
put $\lfrak_0:=\t_0\oplus [\lfrak,\lfrak]$.

Let us  recall some results from~\cite{alg_hamil}, Subsection 5.2.
Put $\a_{G,X}^{(X_L)}=\overline{\im\mu_{Z(L)^\circ,
X_L}}\index{agxxl@$\a_{G,X}^{(X_L)}$}$. The restriction of
$\pi_{L,\lfrak}$ to $\z(\lfrak)=\lfrak^L$ is a closed embedding, so we may
consider $\z(\lfrak)$ as a subvariety of $\lfrak\quo L$. Under this
identification, $\a_{G,X}^{(X_L)}$ is identified with
$\overline{\im\psi_{L,X_L}}$. It turns out that $\a_{G,X}^{(X_L)}$
is an affine subspace in $\z(\lfrak)$ of dimension $\defe_G(X)$
intersecting  $\t_0$ in a unique point  $\xi_0$. Taking $\xi_0$ for
0 in the affine space $\a_{G,X}^{(X_L)}$, we may (and will) consider
$\a_{G,X}^{(X_L)}$ as a vector space.

\begin{defi}\label{Def:2.5.1}
The vector space $\a_{G,X}^{(X_L)}$ is called  the {\it Cartan
space} of   $X$ associated with $X_L$.
\end{defi}

The affine subspace $\a_{G,X}^{(X_L)}\subset \z(\lfrak)$ is stable under
the natural action $N_G(L,X_L):\z(\lfrak)$. The point $\xi_0$ is
$N_G(L,X_L)$-invariant. Denote the image of $N_G(L,X_L)$ in
$\GL(\a_{G,X}^{(X_L)})$ by
$W_{G,X}^{(X_L)}\index{wgxxl@$W_{G,X}^{(X_L)}$}$. We remark that
this definition of $W_{G,X}^{(X_L)}$ differs slightly from that
given in~\cite{alg_hamil}. However, if  $G$ is connected, then the
both definitions coincide and $W_{G,X}^{(X_L)}\cong N_G(L,X_L)/L$.

\begin{defi}\label{Def:2.5.2}
The linear group $W_{G,X}^{(X_L)}$ is said to be the {\it Weyl
group} of the Hamiltonian $G$-\!\! variety $X$ associated with
$X_L$.
\end{defi}

The previous definition in the case of cotangent bundles is due to
Vinberg, \cite{Vinberg2}.

\begin{Rem}\label{Rem:2.5.3}
The affine map $\xi\mapsto \xi-\xi_0:\a_{G,X}^{(X_L)}\rightarrow
\z(\lfrak)$ induces an isomorphism $\a_{G,X}^{(X_L)}\cong \z(\lfrak)\cap
\t_0^\perp\cong \z(\lfrak)/\t_0$. This isomorphism is
$N_G(L,X_L)/L$-equivariant. Note that $\a_{G,X}^{(X_L)}$ is
identified with the Lie algebra of the group $A_{G,X}^{(X_L)}\cong
Z(L)^\circ/T_0$.

Further, note that if $X$ posesses an action of $\K^\times$
satisfying the assumptions of Lemma \ref{Lem:2.3.8}, then $0\in
\a_{G,X}^{(X_L)}$ or, in other words, $\a_{G,X}^{(X_L)}=\z(\lfrak)\cap
\lfrak_0^\perp$.
\end{Rem}


\begin{defi}\label{Def:2.5.5}
The lattice $\X(A_{G,X}^{(X_L)})\subset (\z(\lfrak)/\t_0)^*\cong
\z(\lfrak)/\t_0\cong \a_{G,X}^{(X_L)}$ is called the {\it weight
lattice} of the Hamiltonian $G$-variety  $X$ associated with $X_L$
and is denoted by
$\X_{G,X}^{(X_L)}\index{xgxxl@$\X_{G,X}^{(X_L)}$}$.
\end{defi}

\begin{defi}\label{Def:2.5.6}
The annihilator of $\A_{G,X}^{(X_L)}$ in
$\X(A_{G,X}^{(X_L)})=\X_{G,X}^{(X_L)}$ is called the {\it root
lattice} of the Hamiltonian  $G$-\!\! variety $X$ associated with
$X_L$ and is denoted by
$\Lambda_{G,X}^{(X_L)}\index{zzzl@$\Lambda_{G,X}^{(X_L)}$}$.
\end{defi}

\begin{Rem}\label{Rem:2.5.7}
Now let $L'$ be another principal centralizer and $X_{L'}$ an
$L'$-cross-section of $X$. Then there is an element $g\in G$  such
that $gLg^{-1}=L', gX_L=X_{L'}$.  Such an element $g$ is defined
uniquely up to the right multiplication by an element from
$N_G(L,X_L)$.  Then $g(\a_{G,X}^{(X_L)})=\a_{G,X}^{(X_{L'})},
gW_{G,X}^{(X_L)}g^{-1}=W_{G,X}^{(X_{L'})},
g\X_{G,X}^{(X_L)}=\X_{G,X}^{(X_{L'})}$. 

Further, let $\varphi$ be a Hamiltonian morphism of $X$ satisfying
the assumptions of Definition \ref{Def:3.8} for $L,X_L$. Then
$\varphi$ preserves $X_{L'}$ and acts on it by $gt_\varphi g^{-1}$,
where $t_{\varphi}$ is such as in Definition \ref{Def:3.8}. It
follows that $g\Lambda_{G,X}^{(X_L)}=\Lambda_{G,X}^{(X_{L'})}$.
\end{Rem}

Let us illustrate the definitions given above with two examples.

\begin{Ex}\label{Ex:2.6.2}
Put $G=\SL_2, X=T^*(G/T)$, where $T$ denotes the maximal torus in
$G$ consisting of diagonal matrices. As a $G$-\!\! variety, $X\cong
G*_T(\K e\oplus\K f)$, where
$e=\begin{pmatrix}0&1\\0&0\end{pmatrix},f=\begin{pmatrix}0&0\\1&0\end{pmatrix}$
and  $t\in T\cong \K^\times$ acts on $\K e$ by $t^2$ and on $\K f$
by $t^{-2}$. The moment map is given by $[g, xe+yf]\mapsto
\Ad(g)(xe+yf)$. It is clear that $m_G(X)=\dim G$.
$$\mu_{G,X}^{-1}(\t^{pr})=\{[\begin{pmatrix}a&b\\c&d\end{pmatrix},\begin{pmatrix}0&x\\y&0\end{pmatrix}]|
ad=-bc=\frac{1}{2}, ya^2=xb^2, x\neq 0\}.$$ Therefore
$\mu_{G,X}^{-1}(\t^{pr})$ is irreducible and $W_{G,X}^{(X_T)}=W(\g)$
for a unique cross-section $X_T$ of $X$. The matrix $-E\in T$ acts
on  $X$ trivially. The action of $T/\{\pm E\}$ on $X_T$ is effective
whence  $\X_{G,X}^{(X_T)}=\Lambda(\g)$. The group $T^{W(\g)}$
consists of two elements. Therefore
$\Lambda_{G,X}^{(X_T)}=\Lambda(\g)$ or $2\Lambda(\g)$. Put
$n=\begin{pmatrix}0&1\\-1&0\end{pmatrix},
h=\begin{pmatrix}i&0\\0&-i\end{pmatrix}$. Since $nh=-hn$, we see
that the translation by $h$ is an $N_G(T)$-automorphism of $X_T$. It
is checked directly that $y\mapsto hy:X_T\rightarrow X_T$ is the
restriction of the automorphism $[g,v]\mapsto
[gn^{-1},nv]:X\rightarrow X$. The last automorphism is Hamiltonian,
see  Lemma~\ref{Lem:3.5}. Thus $\Lambda_{G,X}^{(X_T)}=2\Lambda(\g)$.
\end{Ex}

\begin{Ex}\label{Ex:2.6.4}
Let $X=V_1\oplus V_2$, where $V_1,V_2$ are the two-dimensional
irreducible symplectic $\SL_2$-modules. Choose a basis
$e_{i1},e_{i2}$ of $V_i$. For a symplectic form on $V$ we take a
unique skew-symmetric form $\omega$ such that $V_1,V_2$ are
isotropic and
$\omega(xe_{11}+ye_{12},ze_{21}+te_{22})=2\det\begin{pmatrix}
x&z\\y&t\end{pmatrix}$. It is easily seen that $m_G(X)=\dim G$. The
moment map is given by
$$\mu_{G,X}(xe_{11}+ye_{12}+ze_{21}+te_{22})=\begin{pmatrix}xt+yz&2yt\\-2zx&-xt-yz\end{pmatrix}.$$
Therefore
$$\mu_{G,X}^{-1}(\t^{pr})=\{xe_{11}+ye_{12}+ze_{21}+te_{22}|
yt=zx=0, xt+yz\neq 0\}.$$ The variety $\mu_{G,X}^{-1}(\t^{pr})$ has
two connected components $\{xe_{11}+ye_{12}+ze_{21}+te_{22}| y=z=0,
xt\neq 0\}$ и $\{xe_{11}+ye_{12}+ze_{21}+te_{22}| x=t=0, yz\neq
0\}$. Take the first one  for $X_T$. One gets
$W_{G,X}^{(X_T)}=\{1\}, \X_{G,X}^{(X_T)}=\X_G$. The automorphism
$(v_1,v_2)\mapsto (tv_1,t^{-1}v_2), v_1\in V_1,v_2\in V_2,t\in
\K^\times$ is central for any $t\in\K^\times$. Thence
$\Lambda_{G,X}^{(X_T)}=\{0\}$.
\end{Ex}

Our next task is to establish some  properties of
$\a_{G,X}^{(\cdot)},W_{G,X}^{(\cdot)},\Lambda_{G,X}^{(\cdot)},\X_{G,X}^{(\cdot)}$.

\begin{Lem}\label{Lem:2.7.9}
Let $L,X_L$ be such as above. Then
$N_{G^\circ}(L,X_L)=N_G(L,X_L)\cap G^\circ, N_G(L,X_L)G^\circ=G$.
Thus $W_{G^\circ,X}^{(X_L)}$ is a normal subgroup of
$W_{G,X}^{(X_L)}$ and there is a natural epimorphism
$G/G^\circ\rightarrow W_{G,X}^{(X_L)}/W_{G^\circ,X}^{(X_L)}$.
\end{Lem}
\begin{proof}
This follows easily from the observation that $N_{G^\circ}(L,X_L)$
permutes the components of $\mu_{G,X}^{-1}(\lfrak^{pr})$ transitively.
\end{proof}

The following lemma compares the quadruples
$(\a_{G,X}^{(X_L)},W_{G,X}^{(X_L)},\X_{G,X}^{(X_L)},\Lambda_{G,X}^{(X_L)})$
and $(\a_{M,X_M}^{(X_L)},$ $W_{M,X_M}^{(X_L)},
\X_{M,X_M}^{(X_L)},\Lambda_{M,X_M}^{(X_L)})$, where $M$ is a Levi
subgroup in $G$ containing $L$ and $X_M$ is a unique
$M$-cross-section of $X$ containing $X_L$.

\begin{Lem}\label{Lem:2.7.3}
Suppose $G$ is connected. Let $M,X_M$ be such as above. Then
\begin{enumerate}
\item $\a_{M,X_M}^{(X_L)}=\a_{G,X}^{(X_L)},\X_{M,X_M}^{(X_L)}=\X_{G,X}^{(X_L)}, A_{M,X_M}^{(X_L)}=A_{G,X}^{(X_L)}$.
\item $N_M(L,Y)=N_G(L,Y)\cap M$. Equivalently, $W_{M,X_M}^{(X_L)}=
W_{G,X}^{(X_L)}\cap M/L$.
\item  $\Lambda_{M,X_M}^{(X_L)}\subset
\Lambda_{G,X}^{(X_L)}$.
%
%
\end{enumerate}
\end{Lem}
\begin{proof}
The first two assertions follow directly from Definitions
\ref{Def:2.5.1},\ref{Def:2.5.2}.

Proceed to assertion 3. The subvariety $X_M\subset X$ is a unique
component of $\mu_{G,X}^{-1}(\m^{pr})$ containing $X_L$. Thus $X_M$
is $\A_{G,X}^{(\cdot)}$-stable. Let $\varphi\in \A_{G,X}^{(\cdot)}$.
Clearly, $\varphi|_{X_M}$ is a central Hamiltonian automorphism of
$X_M$.  The map $\varphi\mapsto \varphi|_{X_M}$ defines the
homomorphism $\A_{G,X}^{(\cdot)}\rightarrow \A_{M,X_M}^{(\cdot)}$
such that the following diagram is commutative.

\begin{picture}(60,30)
\put(2,20){$\A_{G,X}^{(\cdot)}$}\put(44,20){$\A_{M,X_M}^{(\cdot)}$}
\put(15,2){$A_{G,X}^{(X_L)}\cong A_{M,X_M}^{(X_L)}$}
\put(6,18){\vector(1,-1){11}} \put(46,18){\vector(-1,-1){11}}
\put(11,22){\vector(1,0){32}} \put(5,10){\tiny
$\iota_{G,X}^{(X_L)}$} \put(40,10){\tiny $\iota_{M,X_M}^{(X_L)}$}
\end{picture}

This completes the proof.
\end{proof}

\begin{Lem}\label{Lem:2.7.4}
Let $G$ be connected, $X,X'$ Hamiltonian $G$-varieties and
$\varphi:X\dashrightarrow X'$ a dominant  generically finite
$G$-equivariant rational mapping such that
$\mu_{G,X'}\circ\varphi=\mu_{G,X}$.  Let $X_L'$ be the
$L$-cross-section of  $X'$ described in assertion 2 of
Lemma~\ref{Lem:2.3.7}. Then $\a_{G,X}^{(X_L)}=\a_{G,X'}^{(X_L')},
W_{G,X}^{(X_L)}\subset W_{G,X'}^{(X_L')}.$ If $\varphi$ is
birational, then $\X_{G,X}^{(X_L)}=\X_{G,X'}^{(X_L')},
W_{G,X}^{(X_L)}=W_{G,X'}^{(X_L')}$.
%
\end{Lem}
\begin{proof}
This stems
directly from the definition of $X_L'$. 
\end{proof}

\begin{Lem}\label{Lem:2.7.1}
\begin{enumerate}
\item $\X_{G,X}^{(X_L)}=\X_{G^\circ,X}^{(X_L)}$.
\item The lattices $\X_{G,X}^{(X_L)},\Lambda_{G,X}^{(X_L)}\subset\a_{G,X}^{(X_L)}$ are $W_{G,X}^{(X_L)}$-invariant. 
\item  $w\xi-\xi\in \Lambda_{G,X}^{(X_L)}$ for all $w\in W_{G,X}^{(X_L)},\xi\in \X_{G,X}^{(X_L)}$.
\end{enumerate}
\end{Lem}
\begin{proof}
The first assertion is obvious. The second one stems from Remark
\ref{Rem:2.5.7}. Proceed to the third assertion. An element $t\in
A_{G,X}^{(X_L)}$ is $N_G(L,X_L)$-invariant iff for all $\xi\in
\X_{G,X}^{(X_L)},w\in W_{G,X}^{(X_L)}$ the equality
$\langle\xi,t\rangle=\langle\xi, w t\rangle$ holds, equivalently,
$\langle w^{-1}\xi-\xi,t\rangle=0$. Since $\A_{G,X}^{(X_L)}\subset
(A_{G,X}^{(X_L)})^{N_G(L,X_L)}$, we are done.
\end{proof}

\begin{Lem}\label{Lem:2.7.5}
Suppose $G$ is  connected and $0\in \overline{\im\psi_{G,X}}$. Let
$G^0$ be a connected subgroup of $G$ containing $(G,G)$.  Put
$L^0=L\cap G^0$. Recall (Lemma~\ref{Lem:2.3.7}) that $L^0$ is the
principal centralizer and $X_L$ is an $L^0$-cross-section of the
$G^0$-variety $X$. Then the following assertions hold.
\begin{enumerate}
\item The inclusion
$N_{G^0}(L^0,X_L)\hookrightarrow N_G(L,X_L)$ induces an isomorphism
$W_{G,X}^{(X_L)}\cong W_{G^0,X}^{(X_L)}$.
\item  $\a_{G,X}^{(X_L)}\cap\g^0\subset \a_{G^0,X}^{(X_L)}$;
 $\a_{G,X}^{(X_L)}\cap\g^0$ is a $W_{G,X}^{(X_L)}$-stable subspace of $\a_{G^0,X}^{(X_L)},\a_{G,X}^{(X_L)}$.
\item  $W_{G,X}^{(X_L)}$ acts trivially on
$\a_{G^0,X}^{(X_L)}/(\a_{G,X}^{(X_L)}\cap\g^0)$,
$\a_{G,X}^{(X_L)}/(\a_{G,X}^{(X_L)}\cap\g^0)$.
\item The orthogonal projection  $\g\rightarrow \g^0$ induces
the projection $p:\a_{G,X}^{(X_L)}\rightarrow \a_{G^0,X}^{(X_L)}$,
that is $W_{G,X}^{(X_L)}$-equivariant. 
%
\item $\Lambda_{G,X}^{(X_L)},\Lambda_{G^0,X}^{(X_L)}\subset
\a_{G,X}^{(X_L)}\cap\g^0$, and
$\Lambda_{G,X}^{(X_L)}=\Lambda_{G^0,X}^{(X_L)}$ (the equality of
lattices in $\a_{G,X}^{(X_L)}\cap\g^0$).
\end{enumerate}
\end{Lem}
\begin{proof}
Assertion 1 stems directly from $N_G(L,X_L)\cap
G^0=N_{G^0}(L^0,X_L)$.

To prove the remaining assertions we may replace $G$ by a covering
and assume that $G=Z\times G^0$, where $Z$ is a central subgroup of
$G$ with the Lie algebra $\g^{0\perp}$.

Let $T_0,T_0^0$ denote the inefficiency kernels of the actions
$Z(L)^\circ,Z(L^0)^\circ:X_L$. Clearly, $\t_0^0=\t_0\cap \g^0$. Note
that $\a_{G,X}^{(X_L)}=\z(\lfrak)\cap\t_0^\perp,
\a_{G^0,X}^{(X_L)}=\z(\lfrak^0)\cap \t_0^{0\perp}$, since $0\in
\overline{\im\psi_{G,X}},0\in\overline{\im\psi_{G^0,X}}$ (see
Remark~\ref{Rem:2.5.3}). Assertion 2 stems from the equalities
$\a_{G,X}^{(X_L)}\cap\g^0=\z(\lfrak)\cap\g^0\cap \t_0^\perp,
\a_{G^0,X}^{(X_L)}=\z(\lfrak^0)\cap\t_0^{0\perp}=\z(\lfrak)\cap\g^0\cap
(\t_0\cap\g^0)^\perp$.

By definition,  $W_{G,X}^{(X_L)}$ (resp., $W_{G^0,X}^{(X_L)}$) acts
on $\a_{G,X}^{(X_L)}$ (resp., on $\a_{G^0,X}^{(X_L)}$) as
$N_G(L,X_L)/L$ (resp., as $N_{G^0}(L^0,X_L)/L^0$). But both
$N_G(L,X_L),N_{G^0}(L^0,X_L)$ preserve $\g^0$ and act trivially
$\g/\g^0$.

From the definitions of the corresponding moment maps we see that
$\mu_{Z(L^0)^\circ,X_L}$ is the composition of
$\mu_{Z(L)^\circ,X_L}$ and the orthogonal projection
$\g\twoheadrightarrow \g^0$. Therefore
$\a_{G^0,X}^{(X_L)}=p(\a_{G,X}^{(X_L)})$. Since the orthogonal
projection $\lfrak\rightarrow \lfrak^0$ is $N_G(L,X_L)$-equivariant, we see
that $p$ is $W_{G,X}^{(X_L)}$-equivariant. 

Finally, proceed to assertion 5. At first, we describe the relation
between $\A_{G,X}^{(X_L)}$ and $\A_{G^0,X}^{(X_L)}$. Namely, choose
$l\in Z(L)^\circ=Z\times Z(L^0)^\circ$, $l=(z,l^0)$. Let us show
that the following conditions are equivalent.
\begin{itemize}
\item[(a)] The translation by $l$ in $X_L$ is the restriction of a
$G$-equivariant automorphism of $X$.
\item[(b)] The translation by $l^0$ in $X_L$ is the restriction of a
$G^0$-equivariant automorphism of $X$.
\end{itemize}
To prove $(a)\Rightarrow (b)$ note that the translation by an
element of $Z$ in $X$ is a $G$-\!\! equivariant automorphism.
Conversely, let (b) hold and $\varphi$ be the corresponding
$G^0$-equivariant automorphism. The automorphism $\varphi$ is
Hamiltonian and central, thanks to Lemma~\ref{Lem:3.3}. The group
$Z$ acts on $X$  by $G^0$-equivariant automorphisms leaving $X_L$
stable. Therefore $\varphi$ is $Z$-equivariant (Lemma~\ref{Lem:3.4})
whence $G$-\!\! equivariant.

The lattice $\Lambda_{G,X}^{(X_L)}\subset \X(Z(L)^\circ)$ (resp.
$\Lambda_{G^0,X}^{(X_L)}\subset \X(Z(L^0)^\circ)$) is the
annihilator of the group of all elements $l\in Z(L)^\circ$ (resp.,
$l^0\in Z(L^0)^\circ$) satisfying (a) (resp., (b)). Note that
$\X(Z(L^0)^\circ)$ can be identified with the annihilator of
$Z\subset Z(L)^\circ$ in $\X(Z(L)^\circ)$. The equivalence of
(a),(b) yields $\Lambda_{G,X}^{(X_L)}=\Lambda_{G^0,X}^{(X_L)}$ (the
equality of subgroups in $\X(Z(L)^\circ)$). An element of $ZT_0$
satisfies (a). It follows that $\Lambda_{G,X}^{(X_L)}$ and
$\Lambda_{G^0,X}^{(X_L)}$ considered as lattices in
$\a_{G,X}^{(X_L)}, \a_{G^0,X}^{(X_L)}$, respectively, lie in the
orthogonal complement to $\z+\t_0$ in $\z(\lfrak)$. This orthogonal
complement coincides with $\a_{G,X}^{(X_L)}\cap \g^0$. Since
$\Lambda_{G,X}^{(X_L)}, \Lambda_{G^0,X}^{(X_L)}$ coincide as
subgroups in $\X(Z(L)^\circ)$, they also coincide as lattices in
$\a_{G,X}^{(X_L)}\cap \g^0$.
\end{proof}

Till the end of the subsection  $G$ is connected and $X$ satisfies
the equivalent conditions of Lemma \ref{Lem:2.3.5}. We denote by $T$
a maximal torus of $G$. Let $X_T$ be a $T$-cross-section of   $X$.
Then $\a_{G,X}^{(X_T)}=\t$ and $W_{G,X}^{(X_T)}\subset W(\g)$. Let
$G_1,\ldots, G_k$ be all simple normal subgroups of $G$ so that
$G=Z(G)^\circ G_1\ldots G_k$. Put $T_i=G\cap T_i$. Denote by
$X_{T_i},i=\overline{1,k}$ the $T_i$-cross-sections  described in
assertion 4 of Lemma~\ref{Lem:2.3.7}.

\begin{Lem}\label{Lem:2.7.8}
We preserve the notation introduced in the previous paragraph.
\begin{enumerate}
\item There is the inclusion $W_{G,X}^{(X_T)}\subset \prod_{i=1}^k
W_{G_i,X}^{(X_{T_i})}$ of subgroups in $W(\g)$. The image of
$W_{G,X}^{(X_T)}$ in $\GL(\t_i)$ coincides with
$W_{G_i,X}^{(X_{T_i})}$.
\item Suppose that the linear group  $W_{G,X}^{(X_L)}\in\GL(\t)$ is generated by reflections.
Then $W_{G,X}^{(X_T)}= \prod_{i=1}^k W_{G_i,X}^{(X_{T_i})}$.
\end{enumerate}
\end{Lem}
\begin{proof}
Changing $G$ by a covering, we may assume that
$G=Z(G)^\circ\times\prod_{i=1}^k G_i$. Denote by $G^{(i)}$ the
product of $Z(G)^\circ$ and $G_j, j\neq i$.

Let us show that $N_G(T,X_T)\subset Z(G)^\circ\times \prod_{i=1}^k
N_{G_i}(T_i,X_{T_i})$. Choose $n\in N_G(T,X_T)$. We have to show
that the projection $n_i$ of $n$ to $G_i$ lies in
$N_{G_i}(T_i,X_{T_i})$. Since $X_T\subset X_{T_i}$, we see that
$nX_{T_i}\cap X_{T_i}\neq \varnothing$. Since $X_{T_i}$ is
$G^{(i)}$-stable, we have $n_iX_{T_i}\cap X_{T_i}\neq \varnothing$.
But the variety $\mu_{G_i,X}^{-1}(\t_i^{pr})$ is smooth and
$n_iX_{T_i},X_{T_i}$ are its irreducible components. So $n_i\in
N_{G_i}(T_i,X_{T_i})$.

Now let us show that the projection  $N_G(T,X_T)\rightarrow
N_{G_i}(T_i,X_{T_i})$ is surjective or, equivalently, that for any
$g\in N_{G_i}(T_i,X_{T_i})$ there exists $h\in G^{(i)}$ such that
$gh\in N_G(T,X_T)$. Recall that the subset $G^{(i)}X_T$ is dense in
$X_{T_i}$. Thus for $y_1\in X_T, h_1\in G^{(i)}$ in general position
there exist $y_2\in X_T,h_2\in G^{(i)}$ such that $gh_1y_1=h_2y_2$.
Since  $g\in G_i, h_2\in G^{(i)}$, we have $gh_2=h_2g$. Thus $gh
y_1=y_2$, where $h=h_2^{-1}h_1$. The inclusions
$\mu_{G,X}(y_1),\mu_{G,X}(y_2)\in \t^{pr}$ imply $gh\in N_G(\t)$.
Thence $X_T,gh X_T$ are irreducible components of
$\mu_{G,X}^{-1}(\t^{pr})$ and $X_T\cap ghX_T\neq \varnothing$. It
follows that $gh X_T=X_T$. This completes the proof of assertion 1.

Proceed to assertion 2. Let $s\in W_{G,X}^{(X_L)}$ be a reflection.
Since $\t_i\subset \t, i=\overline{1,k},$ is  $s$-stable, there
exists $j\in \{1,\ldots,k\}$ such that $s$ acts trivially on $\t_i,
i\neq j$. Therefore there is $g\in G_j\cap N_G(T,X_T)$ mapping to
$s$ under the canonical epimorphism $N_G(T,X_T)\twoheadrightarrow
W_{G,X}^{(X_T)}$. In other words, $s\in W_{G_j,X}^{(X_{T_j})}$.
Since $W_{G,X}^{(X_T)}$ is generated by reflections, we see that
$W_{G,X}^{(X_T)}=\prod_{i=1}^k \Gamma_i$, where $\Gamma_i\subset
W_{G_i,X}^{(X_{T_i})}$. Since the projection
$W_{G,X}^{(X_T)}\rightarrow W_{G_i,X}^{(X_{T_i})}$ is surjective, we
get $\Gamma_i=W_{G_i,X}^{(T_i)}$.
\end{proof}

\begin{Lem}\label{Lem:2.7.6}
Let $G$ be connected,  $X$  satisfy the equivalent conditions of
Lemma \ref{Lem:2.3.5}. Suppose that
\begin{itemize}
\item[(a)]
$W_{G,X}^{(\cdot)}$ is generated by reflections. \item[(b)] All
fibers of $\psi_{G,X}$ have the same dimension.
\item[(c)] If $\varphi\in \A_{G,X^{pr}}^{(\cdot)}$ is such that
$\varphi$ and $\varphi^{-1}$ are defined on all divisors contained
in $X\setminus X^{pr}$, then $\varphi$ is a regular automorphism of
$X$.
\end{itemize}
 Then
$\Lambda_{G,X}^{(X_T)}=\bigoplus_{i=1}^k
\Lambda_{G_i,X}^{(X_{T_i})}$.
\end{Lem}
\begin{proof}
By Lemma~\ref{Lem:2.7.5},
$\Lambda_{G,X}^{(X_T)}=\Lambda_{(G,G),X}^{(X_T)}$. So we may (and
will) assume that  $G$ is semisimple. Replacing $G$ with a covering,
we may assume, in addition, that $G=G_1\times G_2\times\ldots\times
G_k$. For $i=1,\ldots k$ put $G^{(i)}=\prod_{j\neq i}G_i$. The proof
of the lemma is in three steps.

{\it Step 1.} Let us check that $\A_{G_i,X}^{(X_{T_i})}=T_i\cap
\A_{G,X}^{(X_T)}$ and

\begin{equation}\label{eq:4.7:1}
\iota_{G,X}^{(X_T)\,-1}(t)=\iota_{G_i,X}^{(X_{T_i})-1}(t), t\in
\A_{G_i,X}^{(X_{T_i})}.\end{equation}

Let $t\in \A_{G_i,X}^{(X_{T_i})}$. Put
$\varphi:=\iota_{G_i,X}^{(X_{T_i})\,-1}(t)$. Note that $t\in
(T_i)^{W_{G,X}^{(X_{T_i})}}$. By Lemma~\ref{Lem:2.7.8}, $t\in
T^{W_{G,X}^{(X_T)}}$. Put
$\widetilde{\varphi}:=\iota_{G,X^{pr}}^{(X_T)\,-1}(t)$. Let us check
that the rational mappings
$\varphi,\widetilde{\varphi}:X\dashrightarrow X$ coincide. The both
rational mappings are $G_i$-equivariant and their restriction to
$G^{(i)}X_T$ coincide with the translation by $t$. Since
$G_iG^{(i)}X_T=GX_T=X^{pr}$, we have $\varphi=\widetilde{\varphi}$.
By assertion 2 of Lemma~\ref{Lem:3.5}, $\widetilde{\varphi}\in
\A_{G,X}^{(\cdot)}$, in other words,  $t\in \A_{G,X}^{(X_T)}$.

Conversely, choose $t\in T_i\cap \A_{G,X}^{(X_T)}$. Put
$\varphi:=\iota_{G,X}^{(X_T)\,-1}(t)$. Since $\varphi$ is
$G$-equivariant, we have $\varphi(gy)=g\varphi(y)=tgy$ for $g\in
G_i,y\in X_T$. Hence $\varphi(x)=tx$ for any $x\in X_{T_i}$. In
particular, $\varphi\in \A_{G_i,X}^{(\cdot)}$ and $t\in
\A_{G_i,X}^{(X_{T_i})}$.

{\it Step 2.} By the previous step,
$\prod_{i=1}^k\A_{G_i,X}^{(X_{T_i})}\subset \A_{G,X}^{(X_T)}$,
equivalently, $\Lambda_{G,X}^{(X_T)}\subset
\bigoplus_{i=1}^k\Lambda_{G_i,X}^{(X_{T_i})}$. The inverse inclusion
will follow if we show that $t:=(t_1,\ldots,t_k)\in
\A_{G,X}^{(X_T)}$ implies $t_i\in \A_{G_i,X}^{(X_{T_i})}$. Thanks to
step 1, it is enough to show that $t_i\in \A_{G,X}^{(X_T)}$. Since
$t_i\in T^{W_{G,X}^{(X_T)}}$, we have $t_i\in
\A_{G,X^{pr}}^{(X_T)}$. Put
$\varphi_i=\iota_{G,X^{pr}}^{(X_T)\,-1}(t_i)$. The elements
$\varphi_i\in \A_{G,X^{pr}}^{(\cdot)}$  mutually commute and
\begin{equation}\label{eq:4.7:2}
\prod_{i=1}^k \varphi_i=\varphi,
\end{equation}
where $\varphi:=\iota_{G,X}^{(X_T)\,-1}(t)$.

{\it Step  3.} Put $X_i^{pr}=G_iX_{T_i}\subset X$. It is an open
affine $G$-stable subvariety of $X$. Further,
\begin{equation}\label{eq:4.7:3}
X\setminus
X_i^{pr}=\psi_{G_i,X}^{-1}\left((\t_i\setminus\t_i^{pr})/W(\g_i)\right)=\psi_{G,X}\left((\t\setminus
(\t_i^{pr}\oplus \bigoplus_{j\neq
i}\t_j))/W(\g)\right),\end{equation}

It follows from Lemma~\ref{Lem:2.7.8} that $t_i\in
T_i^{W_{G_i,X}^{(X_{T_i})}}$. Applying (\ref{eq:4.7:1}) to
$(G,X_i^{pr})$, we get $\varphi_i=\iota_{G_i,X}^{(X_{T_i})-1}(t_i)$.
Therefore $\varphi_i,\varphi_i^{-1}\in \Aut^G(X_i)$. By
(\ref{eq:4.7:3}) and the assumption that  $\widehat{\psi}_{G,X}$ is
equidimensional, the divisors $X\setminus X_i^{pr},X\setminus
X_j^{pr}$ do not have the same irreducible component provided $i\neq
j$. So $\varphi_j^{-1},j\neq i,$ is defined at the components of
$X\setminus X_i$. Since $\varphi_i=\varphi\prod_{j\neq
i}\varphi_j^{-1}$ (see \ref{eq:4.7:2}), $\varphi_i$ satisfies (c).
By Lemma~\ref{Lem:3.4}, $\varphi_i\in \A_{G,X}^{(\cdot)}$, in other
words, $t_i\in \A_{G,X}^{(X_T)}$.
\end{proof}

\begin{Rem}\label{Rem:2.7.10}
Clearly, condition (c) holds whenever $X$ is affine.  In fact, (b)
also holds if $X$ is affine, see, for example, \cite{Weyl}, Theorem
5.1.1.
\end{Rem}

\section{The case of cotangent bundles}\label{SECTION_cotangent}
In this section $G$ is a   reductive group,  $X_0$ is a smooth
irreducible $G$-variety, $X=T^*X_0$. Our objective here is to
interpret invariants of $X$ introduced in the previous section in
terms of $X_0$.

Till the otherwise is indicated, we suppose $G$ is connected.
 Fix a Borel subgroup $B\subset G$ and a maximal torus $T\subset
B$.

\begin{defi}\label{Def:0.1.2}
The lattice $\X_{G,X_0}:=\{\lambda\in \X(T)|
\K(X_0)^{(B)}_\lambda\neq\{0\}\}$ is called the {\it weight} lattice
of the $G$-variety $X_0$. Put
$\a_{G,X_0}=\X_{G,X_0}\otimes_\Z\K\index{agx@$\a_{G,X_0}$}$. We call
the subspace $\a_{G,X_0}\subset \t^*$  the {\it Cartan space} of
$X_0$. The dimension of $\a_{G,X_0}$ is called the {\it rank} of
$X_0$ and is denoted by $\rank_G(X_0)\index{rgx@$\rank_{G}(X_0)$}$.
\end{defi}

Next we will define the Weyl group of $X_0$. To this end we need the
notion of a {\it central $G$-valuation}.

\begin{defi}\label{Def:0.1.3}
By a $G$-valuation of $X_0$ we mean a discrete $\Q$-valued
$G$-invariant valuation of $\K(X_0)$. A $G$-valuation  is called
{\it central} if it vanishes on $\K(X_0)^B$.
\end{defi}

A central valuation $v$ determines the element $\varphi_v\in
\a_{G,X_0}(\Q)^*$ by $\langle
\varphi_v,\lambda\rangle=v(f_\lambda)$, where $\lambda\in
\X_{G,X_0}$, $f_\lambda\in \K(X_0)^{(B)}_\lambda\setminus \{0\}$.
The element $\varphi_v$ is well-defined because $v$ is central.

Recall that we have fixed a $W(\g)$-invariant scalar product on
$\t(\Q)$. This induces the scalar product on $\a_{G,X_0}(\Q)\cong
\a_{G,X_0}(\Q)^*$.

\begin{Thm}\label{Thm:0.1.4}
\begin{enumerate}
\item The map $v\mapsto \varphi_v$ is injective. Its image is a finitely
generated rational convex cone in $\a_{G,X_0}(\Q)^*$ called the {\it
central valuation cone} of $X_0$.
\item The central valuation cone is simplicial (that is, there are linearly
independent vectors $\alpha_1,\ldots,\alpha_s\in \a_{G,X_0}(\Q)$
such that the cone coincides with $\{x|
\langle\alpha_i,x\rangle\geqslant 0, i=\overline{1,s}\}$). Moreover,
the reflections corresponding to its faces generate a finite group.
This group is called the {\rm Weyl group} of $X$  and is denoted by
$W_{G,X_0}\index{wgx@$W_{G,X}$}$.
\item The lattice $\X_{G,X_0}\subset \a_{G,X}(\Q)$ is $W_{G,X_0}$-stable.
\end{enumerate}
\end{Thm}

The proof of the first part of the theorem is relatively easy. It is
obtained (in a greater generality) in \cite{Knop4}, Korollare 3.6,
4.2, 5.2, 6.5. The second assertion is much more complicated. It was
proved by Brion, \cite{Brion}, in the  case when  $X_0$ is
spherical. Knop, \cite{Knop4}, used Brion's result to prove the
assertion in the general case. Later, he gave an alternative proof
in \cite{Knop3}. The third assertion of Theorem \ref{Thm:0.1.4}
follows easily from the construction of the Weyl group in
\cite{Knop3}.

Proceed to the definition of the root lattice of $X_0$. First we
recall the notion of a central automorphism. The following
definition was given in \cite{Knop8}.

\begin{defi}\label{Def:6.1} A $G$-automorphism $\varphi$ of $X_0$ is said to be {\it central}
if for any $\lambda$ there exists $a_{\varphi,\lambda}\in \K^\times$
such that $\varphi(f)=a_{\varphi,\lambda}f$ for any $f\in
\K(X_0)^{(B)}_\lambda$. Central automorphisms of $X_0$ form the
group being denoted by $\A_G(X_0)$.
\end{defi}

It turns out that $\A_G(X_0)$ is not a birational invariant of
$X_0$. However, if $X_0^0$ is an open $G$-subvariety in $X_0$, then
$X_0^0$ is $\A_G(X_0)$-stable and the restriction of elements of
$\A_G(X_0)$ to $X_0^0$ induces an embedding $\A_G(X_0)\rightarrow
\A_G(X^0_0)$ (\cite{Knop8}, Theorem 5.1). It turns out that there is
the maximal group among $\A_G(X_0^0)$ (\cite{Knop8}, Theorem 5.4).
We denote this group by $\A_{G,X_0}$. Clearly, this is a birational
invariant of $X_0$.

Put $A_{G,X_0}=\Hom(\X_{G,X_0},\K^\times)$. We define the map
$\iota_{G,X_0}: \A_{G,X_0}\rightarrow A_{G,X_0}$ by
$\lambda(\iota_{G,X_0}(\varphi)):=a_{\varphi,\lambda}$.  Clearly,
$\iota_{G,X_0}$ is a well-defined group homomorphism.

\begin{Prop}[\cite{Knop8}, Theorem 5.5]\label{Prop:6.2}
The map $\iota_{G,X_0}$ is injective and its image is closed.
\end{Prop}

\begin{defi}\label{Def:6.3}
The sublattice $\Lambda_{G,X}\subset \X_{G,X}=\X(A_{G,X})$
consisting of all elements vanishing on $\iota_{G,X}(\A_{G,X})$ is
called the {\it root lattice} of $X$.
\end{defi}

\begin{Rem}\label{Rem:6.0}
Recall that we have fixed the invariant nondegenerate symmetric form
$(\cdot,\cdot)$ on $G$. So we may consider $\a_{G,X_0}$ as a
subspace of $\g$. We also recall  that it is necessary  to fix a
Borel subgroup and a maximal torus in $G$ to define
$\a_{G,X_0},\X_{G,X_0},W_{G,X_0},\Lambda_{G,X_0}$.
\end{Rem}

The main result of this section is a theorem comparing the
quadruples $(\a_{G,X_0},W_{G,X_0},\X_{G,X_0},$ $\Lambda_{G,X_0})$
and
$(\a_{G,X}^{(\cdot)},W_{G,X}^{(\cdot)},\X_{G,X}^{(\cdot)},\Lambda_{G,X}^{(\cdot)})$
in the case when $X_0$ is quasiaffine. To state this theorem we need
some preliminarily considerations.

Now let $G$ be not necessarily connected.
 Till the end of the section we suppose  $X_0$ is quasiaffine. Put
$L:=Z_{G^\circ}(\a_{G^\circ,X_0}), P:=BL$. Since $X_0$ is
quasiaffine, we see that $\a_{G^\circ,X_0}$ is spanned by dominant
weights. Therefore $P$ is a parabolic subgroup of $G$.

For a parabolic subgroup $Q\subset G$ containing $P$ and an open
$Q$-stable subvariety $X_0'\subset X_0$ put
\begin{equation}\label{eq:3.1:1}\N(Q,X'_0)=\{(x,\beta)\in T^*X_0'| \langle
\beta,\Rad_u(\q)x\rangle\}=\mu_{G,X}^{-1}(\q)\cap
T^*X_0'.\index{ngqx0@$\N(Q,X_0')$}\end{equation}  It follows
directly from (\ref{eq:3.1:1}) that $\N(Q,X_0')\subset X'$ is
$N_G(Q)$-stable.

\begin{Thm}[The comparision theorem]\label{Thm:6.4}
Let $G,X_0,L,P$ be such as above. Then the following assertions
hold.
\begin{enumerate}
\item $L$ is the principal centralizer of $X$.
\item There is an open $P$-subvariety $X_0'\subset X_0$ such that
there exists a unique $L$-cross-section $\Sigma$ of $X$ with
$\Sigma\cap \N(P,X_0')\neq \varnothing$.
\item If $G$ is connected, then $\a_{G,X_0}=\a_{G,X}^{(\Sigma)},A_{G,X_0}=A_{G,X}^{(\Sigma)},
\X_{G,X_0}=\X_{G,X}^{(\Sigma)}, W_{G,X_0}=W_{G,X}^{(\Sigma)},$ $
\Lambda_{G,X_0}=\Lambda_{G,X}^{(\Sigma)}$.
\end{enumerate}
\end{Thm}

\begin{Cor}\label{Cor:6.5}
In the notation of the previous theorem, $W_{G,X}^{(\Sigma)}$ is
generated by reflections.
\end{Cor}

If $X_0''\subset X_0'$ are open $P$-subvarieties of $X_0$, then
$\N(P,X_0'')\subset \N(P,X_0')$. It follows that the
$L$-cross-section $\Sigma$ is uniquely determined.

Actually, only the equality of the root lattices is a new result.
The other claims essentially follow from \cite{Knop3}.

The proof of the theorem will be given after some auxiliary
considerations.

\begin{Prop}\label{Prop:1.1.7}
Suppose $G$ is connected. Then $P$ coincides with the intersection
of the stabilizers of all $B$-divisors on $X$. There is an open
$P$-subvariety $X_0'\subset X_0$ and a closed $L$-subvariety
$S\subset X_0'$ such that the morphism $P*_{L}S\rightarrow X_0',
[p,x]\mapsto px, p\in P, x\in S,$ is an isomorphism. For such a
subvariety $S$ the group $L_0$ is the inefficiency kernel of the
action $L:S$.
\end{Prop}
\begin{proof}
This was proved in \cite{Knop4}, Section 2, Lemma 3.1.
\end{proof}

Let $Q$  be a parabolic subgroup of $G$  containing $P$ and $M$ the
Levi subgroup of $Q$ containing $L$.   Put $\widetilde{M}=
N_{G}(Q)\cap N_G(M)$. Clearly, $\widetilde{M}^\circ=M$,
$N_G(Q)=\widetilde{M}\rightthreetimes \Rad_u(Q)$. In the sequel we
will need to reduce the study of the action $G:X_0$ to the study of
some action of $\widetilde{M}$.

\begin{defi}\label{Def:1.1.9}
We say that a triple $(X_0', Z_0,\pi)$ consisting of an open
$N_G(Q)=N_G(\Rad_u(\q))$-stable subvariety $X_0'\subset X_0$, an
$\widetilde{M}$-variety $Z_0$ and an $N_G(Q)$-morphism
$\pi:X_0'\rightarrow Z_0$ (we consider $Z_0$ as a $N_G(Q)$-variety
assuming that $\Rad_u(Q)$ acts trivially on $Z_0$) {\it reduces} the
action $Q:X_0$ if the following conditions are satisfied:
\begin{itemize}
\item[(a)] $X_0'$ is smooth.
\item[(b)] The action $\Rad_u(Q):X_0'$ is free.
\item[(c)] The morphism $\pi$ is smooth and any its fiber consists of a unique $\Rad_u(Q)$-orbit.
\item[(d)] $Z_0$ is  smooth and quasiaffine.
\end{itemize}
\end{defi}


The following lemma is quite standard.

\begin{Lem}\label{Lem:1.1.11} Let $Q,M$ be such as above. Then there exists a triple $(X_0',Z_0,\pi)$
reducing the action $Q:X_0$.
\end{Lem}
\begin{proof}
Put $X^1_0=\{x\in X_0| \dim \Rad_u(Q)x=\dim \Rad_u(Q)\}$. This is an
open $N_G(Q)$-subvariety in $X_0$. From $P\subset Q$ it follows that
$\Rad_u(Q)\subset \Rad_u(P)$. By Proposition \ref{Prop:1.1.7},  the
action  $\Rad_u(P):X_0$ is generically free. Thus $X_0^1\neq
\varnothing$. The action $\Rad_u(Q):X_0^1$ is free because
$\Rad_u(Q)$ is unipotent.

Since $\Rad_u(Q)$ is a unipotent group,
$\Quot(\K[X_0]^{\Rad_u(Q)})=\K(X_0)^{\Rad_u(Q)}$. There is a
finitely generated subalgebra $A\subset \K[X_0]^{\Rad_u(Q)}$ such
that $\Quot(A)=\K(X_0)^{\Rad_u(Q)}$. Changing $A$ by the subalgebra
generated by $\widetilde{M}A$, we may assume that $A$ is
$\widetilde{M}$-stable. Put $\widetilde{Z}_0:=\Spec(A)$. Let $\pi$
denote the morphism $X_0\rightarrow \widetilde{Z}_0$ induced by the
inclusion $A\subset \K[X_0]$.

There is an open $N_G(Q)$-subvariety $X_0^2\subset X_0^1$  such that
$\pi(X_0^2)\subset \widetilde{Z}_0^{reg}$ and $\pi|_{X_0^2}$ is smooth. 
Since $\K(Z_0)= \K(X_0)^{\Rad_u(Q)}$, a general fiber of
$\pi|_{X_0^{2}}$ consists of a unique $\Rad_u(Q)$-orbit. In other
words, there is an open subvariety $\widehat{Z}_0\subset
\widetilde{Z}_0^{reg}$ such that $\pi|_{X_0^{2}}^{-1}(z)$ is a
$\Rad_u(Q)$-orbit for any $z\in Z_0'$. Changing $\widehat{Z}_0$ by
$\widetilde{M}\widehat{Z}_0$, we may assume that $\widehat{Z}_0$ is
an $\widetilde{M}$-variety. We take
$\pi|_{X_0^{2}}^{-1}(\widehat{Z}_0)$ for $X_0'$. Since $\pi|_{X_0'}$
is a smooth morphism, $Z_0:=\pi(X_0')$ is an open subvariety in
$\widehat{Z}_0$.
\end{proof}

\begin{Prop}\label{Prop:3.1.3}
Let $(X_0',\pi,Z_0)$ be a triple reducing the action $Q:X_0$. Then
\begin{enumerate}
\item The natural projection $\N(Q,X_0')\rightarrow X_0'$ is a vector bundle coinciding
with the inverse image of the bundle $Z:=T^*Z_0\rightarrow Z_0$
under the morphism $\pi$. This observation provides the natural
morphism $\widetilde{\pi}:\N(Q,X_0')\rightarrow Z$.
\item There is a unique $M$-cross-section $\Sigma_M\index{zzzs@$\Sigma_M$}$ of the Hamiltonian $G$-variety $X$
such that $\Sigma_M\cap \N(Q,X_0')\neq\varnothing$. The intersection
$\Sigma_M\cap \N(Q,X_0')$ is dense in $\Sigma_M$ and $\Rad_u(Q)
(\Sigma_M\cap \N(Q,X_0'))$ is dense in $\N(Q,X_0')$. The subvariety
$\Sigma_M\subset X$ is $\widetilde{M}$-stable.
\item The restriction  of $\widetilde{\pi}:\N(Q,X_0')\rightarrow Z$ to $\Sigma_M\cap \N(Q,X_0')$
is dominant, injective, $\widetilde{M}$-equivariant. Moreover,
\begin{equation}\label{eq:moments11}
\mu_{M,Z}\circ\widetilde{\pi}=\mu_{M,\Sigma_M}.
\end{equation}
\end{enumerate}
\end{Prop}
\begin{proof}
Assertion 1 follows easily from assumption (b) of
Definition~\ref{Def:1.1.9}.

Proceed to assertion 2. Since $\p\subset \q$, we see that
$N(Q,X_0')\cap \N(P,X_0'')$ is dense in $\N(P,X_0'')$ for any subset
$X_0''\subset X_0$ that can be included into a triple reducing the
action $P:X_0$. From the proof of Theorem 3.2 from \cite{Knop3} it
follows that there is $y\in \N(Q,X_0')\cap N(P,X_0'')$ such that
$\mu_{G,X}(y)\in \a_{G,X}\cap\lfrak^{pr}$. In particular,
$\N(Q,X_0')\cap\mu_{G,X}^{-1}(\m^{pr})\neq\varnothing$. It is well
known that the morphism $\Rad_u(Q)\times \m^{pr}\rightarrow
\m^{pr}+\Rad_u(\q), (q,\xi)\mapsto q\xi$ is an isomorphism. Thus the
morphism $\Rad_u(Q)\times (\N(Q,X_0')\cap
\mu_{G,X}^{-1}(\m^{pr}))\rightarrow \N(Q,X_0'), (g,x)\mapsto gx,$ is
an open embedding with the image $X'\cap
\mu_{G,X}^{-1}(\m^{pr}+\Rad_u(\q))$. It follows that $\N(Q,X_0')\cap
\mu_{G,X}^{-1}(\m^{pr})$ is an irreducible variety. Therefore there
is a unique $M$-cross-section $\Sigma_M$ of $X$ such that
$\Sigma_M\cap\N(Q,X_0')\neq\varnothing$. It remains to check that
$\Sigma_M$ is $\widetilde{M}$-stable. This stems from the
observation that $g\Sigma_M$ is an $M$-cross-section of $X$ and
$\N(Q,X_0')\cap g\Sigma\neq\varnothing$ whenever $g\in
\widetilde{M}$.

Proceed to assertion  3. The morphism in consideration
 is $\widetilde{M}$-equivariant because so is $\widetilde{\pi}$.
 Let us check (\ref{eq:moments11}). It is enough to show that
 the projection of $\mu_{G,X}(y)\in \q$ to $\m$ coincides
 with $\mu_{M,Z}(\widetilde{\pi}(y))$ for any $y\in \N(Q,X_0')$. Put
$y=(x,\beta)$, where $x\in X_0', \beta\in T^*_xX_0$. The annihilator
of $\Rad_u(\q)_*x$ of  $T^*_xX_0$ is naturally identified with
$T^*_{\pi(x)}Z_0$. By the definition of $\widetilde{\pi}$,
$\widetilde{\pi}(y)=(\pi(x),\beta)$. Let $\xi$ denote the projection
of  $\mu_{G,X}(y)\in\q$ to $\m$. Then $(\xi,\eta)=\langle\eta_*x,
\beta\rangle=\langle
\eta_*\pi(x),\beta\rangle=(\mu_{M,Z}(\widetilde{\pi}(y)),\eta)$ for
all $\eta\in\m$.  (\ref{eq:moments11}) is proved.

Let us check that $\widetilde{\pi}|_{\Sigma_M\cap \N(Q,X_0')}$ is
injective. 
%
%
Recall that the action $\Rad_u(Q):X_0'$ is free and that
$\widetilde{\pi}:\N(Q,X_0')\rightarrow Z$ is fiberwise linear and
induces an isomorphism of fibers. Since any fiber of
$\pi=\widetilde{\pi}|_{X_0'}$ consists of one $\Rad_u(Q)$-orbit, a
fiber of  $\widetilde{\pi}:\N(Q,X_0')\rightarrow Z$ has the same
property. Choose $y_1,y_2\in \Sigma_M\cap \N(Q,X_0')$ such that
$\widetilde{\pi}(y_1)=\widetilde{\pi}(y_2)$. It follows from the
last equality that $\mu_{G,X}(y_1)=\mu_{G,X}(y_2)$ and that there
exists $g\in \Rad_u(Q)$ such that $gy_1=y_2$. In particular,
$\Ad(g)\mu_{G,X}(y_1)=\mu_{G,X}(y_2)$. The latter is impossible,
since $\mu_{G,X}(y_i)\in \m^{pr},i=1,2$.

To prove that  $\widetilde{\pi}|_{\Sigma_M\cap \N(Q,X_0')}$ is
dominant it is enough to show that $\dim \Sigma_M=\dim Z$. But $\dim
Z=2\dim Z_0=2(\dim X_0'-\dim \Rad_u(Q))=\dim X-\dim G+\dim M=\dim
\Sigma_M$.
\end{proof}

\begin{Rem}\label{Rem:3.1.4}
The cross-section $\Sigma_M$ from the previous proposition does not
depend on the choice of $X_0'$. Indeed, if $X_{01}',X_{02}'$ are
subsets of $X_0$ included into triples reducing the action $Q:X$,
then  their intersection, being an open $N_G(Q)$-subvariety of
$X_{01}'$, possesses the same property. Now let $\Sigma_M$ be the
$M$-cross-section of $X$ constructed by $X_{01}'\cap X_{02}'$. Then
$\Sigma_M$ intersects both $\N(Q,X_{01}')$ and $\N(Q,X_{02}')$
nontrivially. In the sequel the $M$-cross-section $\Sigma_M$ will be
called {\it distinguished}. We write $\Sigma$ instead of $\Sigma_L$.
\end{Rem}

One may regard $\X_{G,X}$ as a subgroup in $\X(L)$. Put
$L_0=\bigcap_{\xi\in \X_{G,X}}\ker\chi$.

The following proposition proves all assertions of Theorem
\ref{Thm:6.4} but the equality of the root lattices.

\begin{Prop}\label{Prop:3.1.4}
$L$ is the principal centralizer of the Hamiltonian $G$-variety $X$
and $L_0$ is the inefficiency kernel for the action $L:\Sigma$. If
$G$ is connected, then $\a_{G,X_0}=\a_{G,X}^{(\Sigma)},
\X_{G,X_0}=\X_{G,X}^{(\Sigma)}, W_{G,X_0}=W_{G,X}^{(\Sigma)}$.
Finally, $\Sigma$ is an irreducible component of
$\mu_{G,X}(\a^{pr})$, where $\a^{pr}:=\a_{G,X}\cap\lfrak^{pr}$.
\end{Prop}
\begin{proof}
We may assume that $G$ is connected. Let $X'_0,S$ satisfy the
assumptions of Proposition \ref{Prop:1.1.7} and $\pi:X_0'\cong
P*_LS\cong \Rad_u(P)\times S\twoheadrightarrow S$ be the quotient
morphism for the action $\Rad_u(P):X_0'$. Clearly, the triple
$(X_0',\pi,S)$ reduces the action $P:X_0$.  From
Proposition~\ref{Prop:1.1.7} it follows that $L_{0}$ is the
inefficiency kernel for the action $L:S$ and hence also for the
action  $L:T^*S$. By assertion 3 of Proposition~\ref{Prop:3.1.3},
there is a birational $L$-equivariant rational mapping
$\Sigma_L\dashrightarrow T^*S$. Therefore $L_0$ is the inefficiency
kernel for the action $L:\Sigma_L$. Hence $L$ is the principal
centralizer for the Hamiltonian $L$-variety $\Sigma_L$. Applying
assertion 1 of Lemma~\ref{Lem:2.3.7}, we see that $L$ is the
principal centralizer of $X$. It follows from Remark \ref{Rem:2.5.3}
that $\a_{G,X}^{(\Sigma)}=\lfrak\cap\lfrak_0^\perp=\a_{G,X_0},
\X_{G,X_0}=\X_{G,X}^{(\Sigma)}$. The equality of the Weyl groups was
essentially proved in \cite{Timashev}, Remark 5. The last assertion
is now clear.
\end{proof}

The equality $\Lambda_{G,X_0}=\Lambda_{G,X}^{(\Sigma)}$ will follow
if we prove

\begin{Prop}\label{Prop:3.2.3}
The homomorphism $\varphi\mapsto \varphi_*$ (see Lemma
\ref{Lem:3.3}) induces an isomorphism
$\A_{G,X_0}\cong\A_{G,X}^{(\cdot)}$.  The following diagram is
commutative:

\begin{picture}(60,30)
\put(4,22){$\A_{G,X_0}$}\put(44,22){$\A_{G,X}^{(\cdot)}$}
\put(16,2){$A_{G,X_0}\cong A_{G,X}^{(\Sigma)}$}
\put(12,24){\vector(1,0){30}}\put(22,25){\tiny $\varphi\mapsto
\varphi_*$} \put(7,20){\vector(1,-1){13}}
\put(48,20){\vector(-1,-1){13}} \put(7,13){\tiny $\iota_{G,X_0}$}
\put(44,13){\tiny $\iota_{G,X}^{(\Sigma)}$}
\end{picture}

Moreover, $X$ satisifes condition (c) of Lemma \ref{Lem:2.7.6}.
\end{Prop}

The proof will be given after three auxiliary lemmas.

\begin{Lem}\label{Lem:1.3.2}
The natural embedding $\A_{G}(X_0)\hookrightarrow\A_{G,X_0}$ is an
isomorphism.
\end{Lem}
\begin{proof}
Choose $\gamma\in \A_{G,X_0}$. For any rational $G$-algebra $A$ the
subset  $A^{(B)}\subset \K(X_0)$ is $\gamma$-stable. Thus $A$ is
$\gamma$-stable. Choose an open $G$-embedding $X_0\hookrightarrow
\overline{X}_0$, where $\overline{X}_0$ is an affine $G$-variety,
existing by Theorem 1.6 from \cite{VP}. The subalgebra
$\K[\overline{X}_0]\subset \K(X_0)$ is $\gamma$-stable. Further, any
 $G$-stable ideal $I\subset \K[\overline{X}_0]$ is $\gamma$-stable,
 since
$I=\Span_\K(G I^{(B)})$. Applying this observation to the ideal of
all functions vanishing on $\overline{X}_0\setminus X_0$, we see
that $\gamma\in \A_G(X_0)$.
\end{proof}

\begin{Lem}\label{Lem:1.3.5}
Let $(X_0',\pi,Z_0)$ be a triple reducing the action $P:X_0$ and
$\varphi\in \Aut^G(X_0)$. Suppose $X_0'$ is $\gamma$-stable and
there exists $\varphi_0\in \Aut^L(Z_0)$ such that
$\varphi_0\circ\pi=\pi\circ\varphi$. Then the following conditions
are equivalent.
\begin{enumerate}
\item  $\varphi\in \A_{G,X_0}$.
\item $\varphi_0\in \A_{L,Z_0}$.
\end{enumerate}
Under conditions (1),(2),
$\iota_{G,X_0}(\varphi)=\iota_{L,Z_0}(\varphi_0)$, and
$\varphi_0(z)=\iota_{G,X}(\varphi)z$ for any  $z\in Z_0$.
\end{Lem}
\begin{proof}
Note that $\pi$ induces the  $T$-equivariant group homomorphism
$\pi^*:\K(X_0)^{(B)}\rightarrow \K(Z_0)^{(L\cap B)}$. It follows
that (1)$\Leftrightarrow$(2) and
$\iota_{G,X_0}(\varphi)=\iota_{L,Z_0}(\varphi_0)$. The $L$-varieties
$Z_0$ and $S$ (see Proposition \ref{Prop:1.1.7}) are birationally
equivalent. Therefore $L_0$ is the inefficiency kernel for the
action $L:Z_0$ whence $A_{G,X_0}\cong A_{L,Z_0}\cong L/L_0$. From
the definition of $\iota_{\bullet,\bullet}$ it follows that the
automorphism $(\varphi_0\iota_{L,Z_0}(\varphi_0)^{-1})^*$ of
$\K(Z_0)$ acts trivially on $\K(Z_0)^{(L/L_0)}$. Since the field
$\K(Z_0)$ is generated by $\K(Z_0)^{(L/L_0)}$, we get
$\varphi_0=\iota_{L,Z_0}(\varphi_0)$.
\end{proof}

\begin{Lem}\label{Lem:3.2.4}
Let $T_0$ be a quasitorus (i.e. the direct product of a torus and a
commutative finite group) acting on $X_0$ by $G$-equivariant
automorphisms. Then the following assertion hold:
\begin{enumerate}
\item The structure of a Hamiltonian $G$-variety on $X$ is the
restriction of the structure of a Hamiltonian $G\times T_0$-variety.
\item A triple $(X_0',\pi,Z_0)$ reducing the action $P\times T_0^\circ:X$
also reduces the action $P:X_0$, and $\N(P,X_0')=\N(P\times
T_0,X_0')$.
\item The group $L\times T_0^\circ$ is the principal centralizer
and  $\Sigma$ is the distinguished  $L\times
T_0^\circ$-cross-section of the Hamiltonian $G\times T_0$-variety
$X$.
\item $\Sigma$ is $T_0$-stable and the rational mapping $\Sigma\dashrightarrow T^*Z_0$ from assertion
3 of Proposition \ref{Prop:3.1.3} is $T_0$-\!\! equivariant.
\end{enumerate}
\end{Lem}
\begin{proof}
Assertion 1 follows directly from Example~\ref{Ex:2.1.4}. Put
$G_1=G\times T_0, P_1=P\times T_0^\circ, L_1=L\times T_0^\circ,
\widetilde{L}_1=L\times T_0$. Note that $L_1=L_{G_1^\circ,X_0},
P_1=P_{G_1^\circ,X_0},\widetilde{L}_1=N_{G_1}(P_1)\cap
N_{G_1}(L_1)$. Assertion 2 follows now directly from Definition
\ref{Def:1.1.9}. In the proof of the remaining two assertions one
may assume that  $(X_0',\pi,Z_0)$ reduces the action $P_1:X_0$. In
this case $\N(P,X_0')=\N(P_1,X^0_0)$. Since $\lfrak_1^{pr}=\lfrak^{pr}\times
\t_0$, we see that
$\mu_{G_1,X}^{-1}(\lfrak_1^{pr})=\mu_{G,X}^{-1}(\lfrak^{pr})$. This proves
assertion 3. Assertion 4 stems directly from assertion 3 of
Proposition \ref{Prop:3.1.3}.
\end{proof}

\begin{proof}[Proof of Proposition~\ref{Prop:3.2.3}]
Let us check that $\varphi_*\in \A_{G,X}^{(\cdot)}$ for all
$\varphi\in \A_{G,X_0}$ and that
$\iota_{G,X_0}(\varphi)=\iota_{G,X}^{(\Sigma)}(\varphi_*)$. This
will follow if we check that $\varphi_*(x)=\iota_{G,X}(\varphi)x$
for all $x\in \Sigma$.

 $\A_{G,X_0}$ is a closed subgroup of $A_{G,X_0}$
(Proposition~\ref{Prop:6.2}). Thence $T_0:=\A_{G,X_0}$ is a
quasitorus.  Let  $(X_0',\pi,Z_0)$ be a triple reducing the action
of the parabolic subgroup $P\times T_0^\circ\subset G\times T_0$ on
$X_0$. From Lemma~\ref{Lem:1.3.5} it follows that $\varphi$ acts on
$Z_0$, equivalently, on  $Z:=T^*Z_0$ by the translation by
$\iota_{G,X_0}(\varphi)$. Applying assertion 4 of
Lemma~\ref{Lem:3.2.4}, we get the required claim.

It remains to show that for any $\psi\in \A_{G,X^{pr}}^{(\cdot)}$
such that $\psi,\psi^{-1}$ are defined on all divisors contained in
$X\setminus X^{pr}$ there is $\varphi\in \A_{G,X_0}$ such that
$\psi=\varphi_*$. Recall that there is the natural action
$\K^\times:X$ (see the discussion preceding Lemma \ref{Lem:3.3}).
The rational mapping $\psi:X\dashrightarrow X$ is
$\K^\times$-equivariant. Indeed, by Lemma \ref{Lem:2.3.8}, $\Sigma$
is $\K^\times$-stable. Since $\psi$ is central, we see that its
restriction to $\Sigma$ is $\K^\times$-equivariant. By assertion 3
of Proposition \ref{Prop:2.3.1}, $\overline{G\Sigma}=X$. Thus $\psi$
itself is $\K^\times$-equivariant. Note that $\psi(\K[X])=\K[X]$.
Let $T_0$ denote the subgroup of all $\psi\in
\A_{G,X^{pr}}^{(\cdot)}$ such that $\psi(\K[X])=\K[X]$. By the proof
of Corollary \ref{Cor:3.6}, $T_0$ is a closed subgroup of
$\A_{G,X^{pr}}^{(\cdot)}$. So $T_0$ is a quasitorus. The Lie algebra
of $\A_{G,X^{pr}}^{(\Sigma)}=(A_{G,X_0})^{W_{G,X_0}}$ is
$\a_{G,X_0}^{W_{G,X_0}}$. It follows from \cite{Knop4}, Satz 8.2,
that the Lie algebra $\A_{G,X_0}$ also coincides with
$(\a_{G,X_0})^{W_{G,X_0}}$. By above, any element $\psi\in
T_0^\circ$ has the form $\varphi_*$ for some $\varphi\in
\A_{G,X_0}$.

 Since
$\K[X_0]=\K[X]^{\K^\times}$, we see that $\psi(\K[X_0])=\K[X_0]$ for
any $\psi\in T_0$.
So we have the rational action $T_0:X_0$. Let
$\{\gamma_1,\ldots,\gamma_s\}$ be a finite subset of $T_0$ mapping
surjectively onto $T_0/T_0^\circ$. Denote by $X_0^i$ the domain of
definition of $\gamma_i$ and put $X_0^0=\bigcap_{i=1}^s X_0^i$.
Since $T_0^\circ$ acts on $X_0$ by regular automorphism, we see that
$T_0$ acts on $X_0^0$ by regular automorphisms. Since the action of
$T_0$ commutes with that of $G$, we see that $X_0^0$ is $G$-stable.
By Lemma \ref{Lem:3.3}, $\varphi=(\varphi|_{X_0^0})_*$ for any
$\varphi\in T_0$.
 Let $(X_0',\pi,Z_0)$ be a triple reducing the action of
a parabolic subgroup $P\times T_0^\circ\subset G\times T_0$ on
$X_0^0$. By assertion 4 of Lemma~\ref{Lem:3.2.4}, the group $T_0$
acts on $Z:=T^*Z_0$ or, equivalently, on $Z_0$ by translations by
elements of $L/L_0$. It follows from Lemma~\ref{Lem:1.3.5} that
$\varphi|_{X_0}$ is central.
\end{proof}

\section{Properties of Weyl groups, root and weight lattices
of $G$-varieties}\label{SECTION_properties} In this section $X_0$ is
a smooth quasiaffine $G$-variety, $X=T^*X_0$, $L,L_0,P,\Sigma$ have
the same meaning as in the previous section.

First of all, we note that Theorem \ref{Thm:6.4} allows one to
define $\a_{G,X_0},W_{G,X_0},\X_{G,X_0},\Lambda_{G,X_0}$ even if $G$
is disconnected. Namely, we put $\a_{G,X_0}=\a_{G,X}^{(\Sigma)},
W_{G,X_0}=W_{G,X}^{(\Sigma)}, \X_{G,X_0}=\X_{G,X}^{(\Sigma)},
\Lambda_{G,X_0}=\Lambda_{G,X}^{(\Sigma)}$. Obviously,
$\X_{G,X_0}=\X_{G^\circ,X_0},\a_{G,X_0}=\a_{G^\circ,X_0}$. The
following proposition establishes a relation between
$W_{G,X_0},W_{G^\circ,X}$.

\begin{Prop}\label{Prop:3.4.1}
\begin{enumerate}
\item $N_G(B)\cap N_G(T)\subset N_G(L,\Sigma)$.
In particular, $N_G(B)\cap N_G(L)\subset N_G(\a_{G,X_0})$.
\item Let $\Gamma$ denote the image of
$N_G(B)\cap N_G(L)$ in $\GL(\a_{G,X_0})$. Then
$W_{G^\circ,X_0}\subset W_{G,X_0}$ is a normal subgroup and
$W_{G,X_0}=W_{G^\circ,X_0}\Gamma$.
\item If $\a_{G,X_0}=\t$, then
$W_{G,X_0}=W_{G^\circ,X_0}\leftthreetimes\Gamma$.
\end{enumerate}
\end{Prop}
\begin{proof}
Note that  $\K(X_0)^{(B)}$ is $N_G(B)\cap N_G(T)$-stable whence
$N_G(B)\cap N_G(T)\subset N_G(\a_{G,X_0})$. Therefore $N_G(B)\cap
N_G(T)\subset N_G(L)\cap N_G(B)\subset N_G(P)\cap N_G(L)$. By
assertion 3 of Proposition~\ref{Prop:3.1.3}, $N_G(B)\cap
N_G(T)\subset N_G(L,\Sigma)$. Assertion 2 stems from Lemma
\ref{Lem:2.7.9}. It remains to check that $W_{G^\circ,X_0}\cap
\Gamma=\{1\}$ whenever $\a_{G,X_0}=\t$. Under the last assumption,
$W_{G^\circ,X_0}$ is a subgroup of $W(\g)$. But $W(\g)\cap
\Gamma=\{1\}$ since the dominant Weyl chamber in $\t(\R)$ is
$\Gamma$-stable.
\end{proof}

Till the end of the subsection we suppose that $G$ is connected.
 Let $Q$ be a  parabolic subgroup in $G$ containing
$P$ and $M$ the  Levi subgroup in $Q$ containing $T$. Further, let
$(X_0',\pi,Z_0)$ be a triple reducing the action  $\Rad_u(Q):X_0$.
Being birational invariants, the subspace $\a_{M,Z_0}\subset\t$ and
the group $W_{M,Z_0}$ do not depend on the choice of $Z_0$. We write
$\a_{M,X_0/\Rad_u(Q)},W_{M,X_0/\Rad_u(Q)}$ instead of
$\a_{M,Z_0},W_{M,Z_0}$.

\begin{Prop}\label{Prop:3.4.2}
In the  notation introduced above, $\a_{M,X_0/\Rad_u(Q)}=\a_{G,X_0},
\X_{M,X_0/\Rad_u(Q)}=\X_{G,X_0}$ and
$W_{M,X_0/\Rad_u(Q)}=W_{G,X_0}\cap M/T$.
\end{Prop}
\begin{proof}
The equality of the Cartan spaces and the weight lattices follows
from $\K(X_0)^{(B)}= (\K(X_0)^{\Rad_u(Q)})^{(B\cap M)}$.

 The equality $\a_{M,Z_0}=\a_{G,X_0}$ implies that $P\cap
M=P_{M,Z_0}$ whence there exists a triple $(Z_0',\pi_1,Y_0)$
reducing the action $P\cap M:Z_0$, Lemma~\ref{Lem:1.1.11}. Put
$X_0''=\pi^{-1}(Z_0')$. Let us check that the triple $(X_0'',
\pi_1\circ\pi,Y_0)$ reduces the action
 $P:X$. Clearly,
$X_0''\subset X_0$ is $P$-stable. Since $(\Rad_u(Q))_x=\{1\}$ and
$(\Rad_u(P)\cap M)_{\pi(x)}=\{1\}$ for all $x\in X_0'$, we see that
the action $\Rad_u(P):X_0$ is free.  The variety $Y_0$ and the
morphism $\pi_1\circ\pi:X_0''\rightarrow Y_0$ satisfy condition (c)
of Definition~\ref{Def:1.1.9}.

Let $\Sigma_M$ be the distinguished  $M$-cross-section of  $X$.
Recall that there is the natural embedding
$G*_{N_G(L,\Sigma)}\Sigma\hookrightarrow X$. Comparing the
dimensions we see that $M\Sigma\cong M*_{M\cap N_G(L,\Sigma)}\Sigma$
is a dense subset of some $M$-cross-section of $X$. Since
$\N(P,X_0'')\subset \N(Q,X_0')$, we see that $M\Sigma\cap
\N(Q,X_0')\neq \varnothing$. Therefore  $M\Sigma$ is a dense subset
of $\Sigma_M$. In other words, $\Sigma$ is an $L$-cross-section of
the Hamiltonian  $M$-variety $\Sigma_M$.

Let $\widetilde{\pi}:\N(Q,X_0')\rightarrow Z:=T^*Z_0$ be the natural
morphism defined in assertion 1 of Proposition~\ref{Prop:3.1.3}.
Recall that $\widetilde{\pi}$ induces the birational $M$-equivariant
rational mapping $\Sigma_M\dashrightarrow Z$. It is checked directly
that $\widetilde{\pi}(\N(X_0'',P))=\N(Z_0',P\cap M)$. By assertion 3
of Proposition \ref{Prop:3.1.3}, the subset
$\widetilde{\pi}(\Sigma)$ is dense in the distinguished
$L$-cross-section $\Sigma^Z$ of $Z$. By Lemma~\ref{Lem:2.7.4},
$W_{M,Z}^{(\Sigma^Z)}=W_{M,\Sigma_M}^{(\Sigma)}$. Thanks to
assertion 3 of Lemma~\ref{Lem:2.7.3},
$W_{M,\Sigma_M}^{(\Sigma)}=W_{G,X}^{(\Sigma)}\cap M/L$.
\end{proof}

Now we will give an interpretation of Lemma~\ref{Lem:2.7.5} for
cotangent bundles.
\begin{Prop}\label{Prop:3.4.3}
Let $T_0$ be a torus acting on  $X_0$ by $G$-equivariant
automorphisms. Put $\widetilde{G}=G\times T_0$. Then
$\a_{\widetilde{G},X_0}\cap\g\subset \a_{G,X_0}$ and
$\a_{G,X_0}\cap(\a_{\widetilde{G},X_0}\cap\g)^\perp\subset
\a_{G,X_0}^{W_{G,X_0}}$.
\end{Prop}
\begin{proof}
The assumptions of Lemma~\ref{Lem:2.7.5} hold (with $G$ instead of
$G^0$ and $G\times T_0$ instead of $G$). All assertions of the
proposition follow directly from that lemma.
\end{proof}

Note that the last proposition is a weaker version of results of
Knop, \cite{Knop4}, S\"{a}tze 8.1,8.2.

Our next objective is to reduce the computation of $\X_{G,X_0},
W_{G,X_0}$ to the case when $\a_{G,X_0}=\t$. To this end we need the
notion of the {\it distinguished} component of $X_0^{L_0^\circ}$.

\begin{Prop}\label{Prop:1.5.6}
We use the notation and the conventions of Proposition
\ref{Prop:1.1.7}. Let $L_1$ be a normal subgroup of $L_{0}$. Then
there exists a unique irreducible component $\underline{X}_0\subset
X^{L_1}_0$ such that $\overline{U\underline{X}_0}=X_0$. This
component coincides with $\overline{FS}$, where $F=\Rad_u(P)^{L_1}$.
\end{Prop}
\begin{proof}
 Let $Z$ be a component of
$X_0^{L_1}$. If $Z\cap PS=\varnothing$, then $\overline{UZ}\subset
X_0\setminus PS$. Note that $S\subset X^{L_1}$. Consider a component
$Z$ of $X_0^{L_1}$ containing $S$. Then
$X_0=\overline{P*_LS}=\overline{\Rad_u(P)S}\subset\overline{US}\subset
\overline{UZ}$. Therefore it remains to prove that there exists a
unique component of $X_0^{L_1}$ intersecting $PS$ and that this
component coincides with $\overline{FS}$. The $L_1$-variety $PS\cong
P*_LS$ is isomorphic to $\Rad_u(P)\times S$, where  $L_1$ acts
trivially on $S$ and  by the conjugation on $\Rad_u(P)$. It follows
that $(PS)^{L_1}=FS$. The group $F$ is irreducible because
unipotent. Being an open subset of an irreducible variety, $P*_LS$
is also irreducible. It follows that $S$ is irreducible. So $FS$ is
irreducible.
\end{proof}

\begin{defi}\label{Def:1.5.7}
The component $\underline{X}_0\subset X_0^{L_1}$ satisfying the
assumptions of Proposition \ref{Prop:1.5.6} is said to be {\it
distinguished}.
\end{defi}

The distinguished component for $L_1=L_{0}$ was first considered by
Panyushev in \cite{Panyushev3}.

Let $\underline{X}_0$ be the distinguished component of
$X_0^{L_0^\circ}$. Clearly, $\underline{X}_0$ is quasiaffine. Thanks
to the following lemma, $\underline{X}_0$ is  smooth.

\begin{Lem}\label{Lem:3.5.1}
Let $Y$ be a quasiaffine variety acted on by a reductive group $H$
and $y\in Y^H\cap Y^{reg}$. Then $Y^H$ is smooth in $y$ and
$T_y(Y^H)=(T_yY)^H$.
\end{Lem}
\begin{proof}
Replacing $Y$ with an affine $H$-variety containing $Y$ as an open
$H$-subvariety, we may assume that $Y$ is affine. In this case one
uses Corollary from \cite{VP}, Subsection 6.5.
\end{proof}

Put $\underline{G}:=N_G(L_0^\circ,\underline{X}_0)/L_0^\circ$. It is
a reductive group acting on $\underline{X}_0$. Note that the Lie
algebra $\underline{\g}$ of $\underline{G}$ can be naturally
identified with $\lfrak_0^\perp\cap \n_\g(\lfrak_0)=\g^{L_0^\circ}\cap
\z(\lfrak_0)^\perp\subset \g$. Further, note that $\t\subset
\n_\g(\underline{\g})$. For a Borel (resp., Cartan) subalgebra in
$\underline{\g}$ we take $\b\cap\underline{\g},
\t\cap\underline{\g}$. We define the triple
$(\a_{\underline{G},\underline{X}_0},
\X_{\underline{G},\underline{X}_0},
W_{\underline{G},\underline{X}_0})$ according to this choice of
Borel and Cartan subalgebras.

\begin{Thm}\label{Thm:3.0.3}
In the above notation,
$\a_{G,X_0}=\a_{\underline{G},\underline{X}_0}=\underline{\t}$,
$W_{G,X_0}=W_{\underline{G},\underline{X}_0}$,
$\X_{G,X_0}=\X_{\underline{G}^\circ,\underline{X}_0}$. 
\end{Thm}
\begin{proof}
Recall that there exist an open $P$-subvariety $X_0'\subset X_0$ and
a closed $L$-subvariety  $S\subset X_0'$ such that $X_0'\cong P*_LS$
(Proposition~\ref{Prop:1.1.7}). Let $\Sigma$ denote the
distinguished  $L$-cross-section of a Hamiltonian $G$-variety $X$.

Note that $\underline{\b}=\p\cap\underline{\g}$.  By Proposition
\ref{Prop:1.5.6}, $S\subset \underline{X}_0$. Consider the morphism
$\varphi: \underline{B}*_{L/L_0^\circ}S \rightarrow X_0'\cap
\underline{X}_0, [g,x]\mapsto gx, g\in \underline{B},x\in S$. Recall
that $X_0'\cong \Rad_u(P)\times S$. Thus $\underline{X}_0\cap
X_0'=\Rad_u(P)^{L_0^\circ} S\cong \Rad_u(P)^{L_0^\circ}\times S$. It
follows that  $\varphi$ is surjective. Since
$\underline{B}=\Rad_u(P)^{L_0^\circ}\leftthreetimes (L/L_0^\circ)$,
we see that $\varphi$ is injective. Being a bijective morphism of
smooth varieties,
$\varphi$ is an isomorphism.  
Thanks to the natural isomorphism $\K(S)^{(B\cap L)}\cong
\K(X_0)^{(B)}$, we get
$\X_{G,X_0}=\X_{L,S}=\X_{\underline{G}^\circ,\underline{X}_0},
\a_{G,X_0}=\a_{L,S}=\lfrak\cap\lfrak_0^\perp=\underline{\t}$. Also we see
that $\underline{X}_0\cap X_0'$ can be included into a triple
reducing the action $\underline{B}:\underline{X}_0$. 

Let $p$ denote the natural projection  $X\rightarrow X_0$. For $x\in
X_0^{L_0^\circ}$ the subspace $p^{-1}(x)^{L_0^\circ}\subset
T^*_xX_0$ is naturally identified with $T_x^*(X_0^{L_0^\circ})$.
Therefore the vector bundle $(T^*X_0)^{L_0^\circ}\rightarrow
X_0^{L_0^\circ}$ coincides with the cotangent bundle of
$X_0^{L_0^\circ}$. In particular, there is the natural embedding
$\underline{X}:=T^*\underline{X}_0\hookrightarrow X$. Let us check
that the following diagram is commutative.
\begin{equation}\label{eq:3.5:2}
\begin{picture}(60,30)
\put(3,2){$\underline{\g}$}\put(2,22){$\underline{X}$}
\put(33,2){$\g$} \put(32,22){$X$} \put(4,20){\vector(0,-1){14}}
\put(34,20){\vector(0,-1){14}} \put(6,4){\vector(1,0){24}}
\put(11,23){\vector(1,0){20}} \put(15,4){\tiny $\hookrightarrow$}
\put(15,23){\tiny $\hookrightarrow$} \put(5,14){\tiny
$\mu_{\underline{G},\underline{X}}$} \put(35,14){\tiny $\mu_{G,X}$}
\end{picture}
\end{equation}

Since $\mu_{G,X}$ is $G$-equivariant, we see that
$\mu_{G,X}(\underline{X})\subset \g^{L_0^\circ}$. By the definition
of $\underline{X}_0$, the group $L_0^\circ$ acts trivially on
$\underline{X}$. Since $\underline{X}_0\subset \underline{X}\cap
\mu_{G,X}^{-1}(0)$, we see that
$\mu_{G,X}(\underline{X})\subset\lfrak_0^\perp$. It follows that
$\mu_{G,X}(\underline{X})\subset \underline{\g}$. The commutativity
of the diagram now stems directly from the definitions of the moment
maps for the actions $\underline{G}:\underline{X},G:X$
(Example~\ref{Ex:2.1.4}).

 Since
$\underline{\b}=\underline{\g}\cap \p$,  (\ref{eq:3.5:2}) yields
\begin{equation}\label{eq:3.5.3}\N(\underline{B},\underline{X}_0\cap X_0')\subset \N(P,X_0').\end{equation}
From the equality
$\a_{\underline{G},\underline{X}_0}=\underline{\t}$ it follows that
$L/L_0^\circ$ is the principal centralizer of the Hamiltonian
$\underline{G}$-variety $\underline{X}$.  By Proposition
\ref{Prop:3.1.4},
\begin{equation}\label{eq:3.5:5} \Sigma\subset X^{L_0^\circ}.
\end{equation} Put $\a^{pr}:=\a_{G,X_0}\cap\lfrak^{pr}$. By Proposition \ref{Prop:3.1.4}, $\Sigma$ is
an irreducible component of  $\mu_{G,X}^{-1}(\a^{pr})$ and the
distinguished $L/L_0^\circ$-cross-section $\underline{\Sigma}$ of
$\underline{X}$ is an irreducible component of
$\mu_{\underline{G},\underline{X}}^{-1}(\a^{pr})$. This observation
and (\ref{eq:3.5.3}) yield $\underline{\Sigma}\subset \Sigma$. From
 $\underline{\Sigma}\subset \Sigma$, (\ref{eq:3.5:5}) and the smoothness of $X^{L_0^\circ}$
 we deduce that $\Sigma\subset \underline{X}$. It follows from (\ref{eq:3.5:2})
 that
$\Sigma\subset \mu_{\underline{G},\underline{X}}^{-1}(\a^{pr})$.
Hence $\underline{\Sigma}=\Sigma$.

Let us check the equality of the Weyl groups. We have to prove that
\begin{equation}\label{eq:3.5:4}
N_{\underline{G}}(L/L_0^\circ,\Sigma)=
N_G(L,\Sigma)/L_0^\circ.\end{equation} Since $L_0$ is the
inefficiency kernel of the action $L:\Sigma$, we get
$N_G(L,\Sigma)\subset N_G(L_0^\circ)$.  Being a unique component of
$X^{L_0^\circ}$ containing $\Sigma$, $\underline{X}$ is
$N_G(L,\Sigma)$-stable. Therefore $N_G(L,\Sigma)\subset
N_G(L_0^\circ,\underline{X})=N_G(L_0^\circ,\underline{X}_0)$.
Clearly,
\begin{equation}\label{eq:3.5:6} N_{\underline{G}}(L/L_0^\circ,\Sigma)=
N_{N_G(L_0^\circ,\underline{X}_0)}(L,\Sigma)/L_0^\circ.\end{equation}
The inclusion $N_G(L,\Sigma)\subset N_G(L_0^\circ,\underline{X}_0)$
and (\ref{eq:3.5:6}) yield (\ref{eq:3.5:4}).
\end{proof}

Now let $G_1,\ldots, G_k$ be all simple normal subgroups in $G$ so
that $G=Z(G)^\circ G_1\ldots G_k$. We take $B_i:=B\cap G_i,
T_i:=T\cap G_i$ for a Borel subgroup and a maximal torus of $G_i$.

\begin{Prop}\label{Prop:3.4.4}
Suppose $\a_{G,X_0}=\t$. Then $\a_{G_i,X_0}=\t_i$,
$W_{G,X_0}=\prod_{i=1}^k W_{G_i,X_0}$, and
$\Lambda_{G,X_0}=\bigoplus_{i=1}^k\Lambda_{G_i,X_0}$.
\end{Prop}
\begin{proof}
Replacing $G$ with a covering, we may assume that
$G=Z(G)^\circ\times G_1\times\ldots\times G_k$. For
$i=\overline{0,k}$ put $B_i=G_i\cap B,
G^{(i)}=Z(G)^\circ\times\prod_{j\neq i} G_j, B^{(i)}=B\cap
G^{(i)},T^{(i)}=T\cap G^{(i)}$.

It follows from Proposition~\ref{Prop:3.1.4} and assertion 4 of
Lemma~\ref{Lem:2.3.7} that $\a_{G_i,X}=\t_i$.

Thanks to Proposition~\ref{Prop:1.1.7},  there is an open
$B$-subvariety $X_0'\subset X_0 $ and a closed $T$-subvariety
$S\subset X_0'$ such that $X_0'\cong B*_TS$. For $i=1,\ldots,k$ put
$Z_{0i}=B^{(i)}S$.  Since $X_0'\cong B*_TS$, we see that
$Z_{0i}\cong B^{(i)}*_{T^{(i)}}S$, $X_0'\cong B_i*_{T_i}Z_{0i}$.
Therefore the subset $X_0'\subset X_0$ can be included into a triple
reducing the action $B_i:X$.

It follows immediately from (\ref{eq:3.1:1}) that
\begin{equation}\label{eq:3.1:3}\mathcal{N}(B,X_0')=\bigcap_{i=1}^k
\mathcal{N}(B_i,X_0').\end{equation}

Let $\Sigma$ be the distinguished $T$-cross-section of $X$ and
$\Sigma_i$ a unique  $T_i$-cross-section of $X$ containing
$G^{(i)}\Sigma$, see assertion 4 of Lemma~\ref{Lem:2.3.7}. Therefore
the proposition will follow if we check the equality $\Sigma_i\cap
\N(B_i,X_0')\neq \varnothing$ and the assumptions of Lemmas
\ref{Lem:2.7.8}, \ref{Lem:2.7.6}. The former inequality stems
directly from $\Sigma\cap \N(B,X^0)\neq \varnothing$ and
(\ref{eq:3.1:3}). By Corollary \ref{Cor:6.5}, $W_{G,X}^{(\Sigma)}$
is generated by reflections. The morphism $\psi_{G,X}:X\rightarrow
\g\quo G$ is equidimensional, thanks to Satz 6.6 from \cite{Knop1}.
Finally, condition (c) of Lemma \ref{Lem:2.7.6} follows from
Proposition \ref{Prop:3.2.3}.
\end{proof}

\bigskip
{\Small Chair of Higher Algebra, Department of Mechanics and
Mathematics, Moscow State University.

\noindent E-mail address: ivanlosev@yandex.ru}
\end{document}